\documentclass[a4paper, reqno,10pt]{amsart}
\usepackage[usenames,dvipsnames]{color}

\usepackage{amsthm,amsfonts,amssymb,amsmath,amsxtra,amsrefs}
\usepackage[all]{xy}
\SelectTips{cm}{}
\usepackage{xr-hyper}
\usepackage[colorlinks=
   citecolor=Black,
   linkcolor=Red,
   urlcolor=Blue]{hyperref}
\usepackage{verbatim}

\usepackage[margin=1.25in]{geometry}

\usepackage{mathrsfs}

\RequirePackage{xspace}
\RequirePackage{etoolbox}
\RequirePackage{varwidth}
\RequirePackage{enumitem}
\RequirePackage{tensor}
\RequirePackage{mathtools}
\RequirePackage{longtable}
\RequirePackage{multirow}

\usepackage[refpage]{nomencl}
\usepackage{ifthen}
\usepackage{lipsum}

\newcommand{\sC}{\ensuremath{\mathscr{C}}\xspace}

\newcommand{\sH}{\ensuremath{\mathscr{H}}\xspace}

\newcommand{\sL}{\ensuremath{\mathscr{L}}\xspace}

\newcommand{\sP}{\ensuremath{\mathscr{P}}\xspace}

\newcommand{\fkm}{\ensuremath{\mathfrak{m}}\xspace}

\newcommand{\BA}{\ensuremath{\mathbb {A}}\xspace}

\newcommand{\BC}{\ensuremath{\mathbb {C}}\xspace}

\newcommand{\BF}{\ensuremath{\mathbb {F}}\xspace}
\newcommand{\BG}{\ensuremath{\mathbb {G}}\xspace}

\newcommand{\BI}{\ensuremath{\mathbb {I}}\xspace}
\newcommand{\BJ}{\ensuremath{\mathbb {J}}\xspace}

\newcommand{\BL}{\ensuremath{\mathbb {L}}\xspace}

\newcommand{\BQ}{\ensuremath{\mathbb {Q}}\xspace}
\newcommand{\BR}{\ensuremath{\mathbb {R}}\xspace}
\newcommand{\BS}{\ensuremath{\mathbb {S}}\xspace}

\newcommand{\BZ}{\ensuremath{\mathbb {Z}}\xspace}

\newcommand{\CA}{\ensuremath{\mathcal {A}}\xspace}

\newcommand{\CE}{\ensuremath{\mathcal {E}}\xspace}

\newcommand{\CH}{\ensuremath{\mathcal {H}}\xspace}

\newcommand{\CL}{\ensuremath{\mathcal {L}}\xspace}
\newcommand{\CM}{\ensuremath{\mathcal {M}}\xspace}
\newcommand{\CN}{\ensuremath{\mathcal {N}}\xspace}
\newcommand{\CO}{\ensuremath{\mathcal {O}}\xspace}

\newcommand{\CX}{\ensuremath{\mathcal {X}}\xspace}

\newcommand{\RG}{\ensuremath{\mathrm {G}}\xspace}
\newcommand{\RH}{\ensuremath{\mathrm {H}}\xspace}

\newcommand{\RN}{\ensuremath{\mathrm {N}}\xspace}

\newcommand{\RR}{\ensuremath{\mathrm {R}}\xspace}

\newcommand{\RT}{\ensuremath{\mathrm {T}}\xspace}
\newcommand{\RU}{\ensuremath{\mathrm {U}}\xspace}

\newcommand{\RZ}{\ensuremath{\mathrm {Z}}\xspace}

\newcommand{\Ad}{{\mathrm{Ad}}}

\newcommand{\Ch}{{\mathrm{Ch}}}

\newcommand{\cl}{{\mathrm{cl}}}

\newcommand{\corr}{\mathrm{corr}}

\newcommand{\del}{\operatorname{\partial Orb}}

\DeclareMathOperator{\diag}{diag}

\DeclareMathOperator{\Eis}{Eis}
\DeclareMathOperator{\End}{End}

\DeclareMathOperator{\Fr}{Fr}

\newcommand{\GL}{\mathrm{GL}}

\newcommand{\GU}{\mathrm{GU}}

\DeclareMathOperator{\Hom}{Hom}

\newcommand{\id}{\ensuremath{\mathrm{id}}\xspace}

\newcommand{\Int}{\ensuremath{\mathrm{Int}}\xspace}

\DeclareMathOperator{\Nm}{Nm}

\DeclareMathOperator{\Orb}{Orb}
\DeclareMathOperator{\ord}{ord}

\newcommand{\PGL}{{\mathrm{PGL}}}
\DeclareMathOperator{\Pic}{Pic}

\renewcommand{\Re}{{\mathrm{Re}}}

\newcommand{\rs}{\ensuremath{\mathrm{rs}}\xspace}

\DeclareMathOperator{\Spec}{Spec}
\DeclareMathOperator{\Spf}{Spf}
\newcommand{\SO}{{\mathrm{SO}}}

\DeclareMathOperator{\sgn}{sgn}

\newcommand{\U}{\mathrm{U}}

\DeclareMathOperator{\vol}{vol}

\newcommand{\Bun}{{\mathrm{Bun}}}

\newcommand{\Sh}{{\mathrm{Sh}}}

\newcommand{\Sht}{{\mathrm{Sht}}}

\newcommand{\Hk}{{\mathrm{Hk}}}

\newcommand{\Mat}{\mathrm{Mat}}

\newcommand{\wt}{\widetilde}
\newcommand{\wh}{\widehat}

\newcommand{\pair}[1]{\langle {#1} \rangle}

\newcommand{\sfrac}[2]{\left( \frac {#1}{#2}\right)}

\newcommand{\ov}{\overline}

\newcommand{\incl}{\hookrightarrow}
\newcommand{\lra}{\longrightarrow}
\newcommand{\imp}{\Longrightarrow}

\newcommand{\bs}{\backslash}

\def\calE{\mathcal{E}}

\renewcommand{\l}{\lambda}

\newcommand{\kbar}{\overline{k}}

\newcommand{\cohoc}[2]{\textup{H}_{c}^{#1}({#2})}     % compact support
\newtheorem{theorem}{Theorem}

\newtheorem{conj}[theorem]{Conjecture}

\newtheorem{thm}[theorem]{Theorem}

\theoremstyle{definition}
\newtheorem{defn}[theorem]{Definition}

\newtheorem{remark}[theorem]{Remark}

\newenvironment{altenumerate}
   {\begin{list}
      {\textup{(\theenumi)} }
      {\usecounter{enumi}
       \setlength{\labelwidth}{0pt}
       \setlength{\labelsep}{0pt}
       \setlength{\leftmargin}{0pt}
       \setlength{\itemsep}{\the\smallskipamount}
       \renewcommand{\theenumi}{\roman{enumi}}
      }}
   {\end{list}}
\newenvironment{altitemize}
   {\begin{list}
      {$\bullet$}
      {\setlength{\labelwidth}{0pt}
	   \setlength{\itemindent}{5pt}
       \setlength{\labelsep}{5pt}
       \setlength{\leftmargin}{0pt}
       \setlength{\itemsep}{\the\smallskipamount}
      }}
   {\end{list}}

\numberwithin{equation}{section}
\numberwithin{theorem}{section}

\newcommand{\sform}{\ensuremath{(\text{~,~})}\xspace}

\setitemize[0]{leftmargin=*,itemsep=\the\smallskipamount}
\setenumerate[0]{leftmargin=*,itemsep=\the\smallskipamount}

\renewcommand{\to}{%
   \ifbool{@display}{\longrightarrow}{\rightarrow}%
   }
\let\shortmapsto\mapsto
\renewcommand{\mapsto}{%
   \ifbool{@display}{\longmapsto}{\shortmapsto}%
   }
\newlength{\olen}
\newlength{\ulen}
\newlength{\xlen}
\newcommand{\xra}[2][]{%
   \ifbool{@display}%
      {\settowidth{\olen}{$\overset{#2}{\longrightarrow}$}%
       \settowidth{\ulen}{$\underset{#1}{\longrightarrow}$}%
       \settowidth{\xlen}{$\xrightarrow[#1]{#2}$}%
       \ifdimgreater{\olen}{\xlen}%
          {\underset{#1}{\overset{#2}{\longrightarrow}}}%
          {\ifdimgreater{\ulen}{\xlen}%
             {\underset{#1}{\overset{#2}{\longrightarrow}}}
             {\xrightarrow[#1]{#2}}}}%
      {\xrightarrow[#1]{#2}}
   }
\makeatother
\newcommand{\xyra}[2][]{%
   \settowidth{\xlen}{$\xrightarrow[#1]{#2}$}%
   \ifbool{@display}%
      {\settowidth{\olen}{$\overset{#2}{\longrightarrow}$}%
       \settowidth{\ulen}{$\underset{#1}{\longrightarrow}$}%
       \ifdimgreater{\olen}{\xlen}%
          {\mathrel{\xymatrix@M=.12ex@C=3.2ex{\ar[r]^-{#2}_-{#1} &}}}%
          {\ifdimgreater{\ulen}{\xlen}%
             {\mathrel{\xymatrix@M=.12ex@C=3.2ex{\ar[r]^-{#2}_-{#1} &}}}
             {\mathrel{\xymatrix@M=.12ex@C=\the\xlen{\ar[r]^-{#2}_-{#1} &}}}}}%
      {\mathrel{\xymatrix@M=.12ex@C=\the\xlen{\ar[r]^-{#2}_-{#1} &}}}%
   }
\makeatletter
\newcommand{\xla}[2][]{%
   \ifbool{@display}%
      {\settowidth{\olen}{$\overset{#2}{\longleftarrow}$}%
       \settowidth{\ulen}{$\underset{#1}{\longleftarrow}$}%
       \settowidth{\xlen}{$\xleftarrow[#1]{#2}$}%
       \ifdimgreater{\olen}{\xlen}%
          {\underset{#1}{\overset{#2}{\longleftarrow}}}%
          {\ifdimgreater{\ulen}{\xlen}%
             {\underset{#1}{\overset{#2}{\longleftarrow}}}
             {\xleftarrow[#1]{#2}}}}%
      {\xleftarrow[#1]{#2}}
   }
\newcommand{\isoarrow}{%
   \ifbool{@display}{\overset{\sim}{\longrightarrow}}{\xrightarrow\sim}%
   }
\renewcommand{\lra}{%
   \ifbool{@display}{\longleftrightarrow}{\leftrightarrow}%
   }

\begin{document}

\thanks{Research of W. Zhang partially supported by the NSF grant DMS $\#$1601144,  and a Simmons fellowship. }

\title[Periods, cycles, and $L$-functions]{
Periods, cycles, and $L$-functions: a relative trace formula approach} 

\author{W. Zhang}
\address{Massachusetts Institute of Technology, Department of Mathematics, 77 Massachusetts Avenue, Cambridge, MA 02139, USA}
\email{weizhang@mit.edu}

\subjclass[2010]{Primary 11F67; Secondary 11F70, 11G40, 14C25, 14G35, 14H60}
\keywords{$L$-functions; Shimura varieties; Drinfeld Shtukas; Special cycles; Relative trace formula; Gross--Zagier formula; Waldspurger formula; Rapoport--Zink spaces; Spherical subgroup}

\date{\today}

\begin{abstract}
This is a report for the author's talk in ICM-2018. 
Motivated by the formulas of  Gross--Zagier and Waldspurger, we review conjectures and theorems on automorphic period integrals, special cycles on Shimura varieties, and their connection to central values of L-functions  and their derivatives. We focus on the global Gan--Gross--Prasad conjectures, their arithmetic versions and some variants in the author's joint work with Rapoport and Smithling. We discuss the approach of relative trace formulas and the arithmetic fundamental lemma conjecture. In the function field setting, Z. Yun and the author obtain a formula for higher order derivatives of L-functions in terms of special cycles on the moduli space of Drinfeld Shtukas.
  
\end{abstract}

\maketitle

\tableofcontents

\section{Introduction}

 We begin with a special case of Pell's equation 
 $$x^2-py^2=1,\quad x,y\in\BZ,
 $$ 
where $p\equiv 1\!\!\!\mod 4$ is a prime number. Let $K$ be the real quadratic field $\BQ[\sqrt{p}]$ and $O_K$ its ring of integers. Then most solutions to Pell's equation can be extracted from the group of units in $O_K$, which is known to be the product of $\{\pm 1\}$ and an infinite cyclic group generated by a fundamental unit. In 1830s, Dirichlet systematically constructed units in $O_K$ using special values of trigonometric functions:
\begin{equation*}\label{cyc unit}
\theta_p=\frac{\prod_{ a\not\equiv\square\!\!\!\!\mod p}\sin\frac{a\pi }{p}}{\prod_{b\equiv\square\!\!\!\!\mod p}\sin\frac{b\pi }{p}},\quad 0<a,b<p/2,
\end{equation*}
where $\square\!\!\mod p$ denotes a square residue. Dirichlet also showed that the obstruction for $\theta_p$ to be a fundamental unit is the class group of $K$, with the help of an infinite series
$$
L\left(s,\sfrac{\cdot}{p}\right)=\sum_{n\geq 1, \, p\nmid n}\sfrac{n}{p} n^{-s},\quad s\in\BC, \quad\Re(s)>1.
$$
Here $\sfrac{\cdot}{p}$ denotes the Legendre symbol for quadratic residues. This is now called a Dirichlet L-series, and it has holomorphic continuation to $s\in\BC$ with a simple zero at $s=0$. What Dirichlet discovered can be stated as two formulas for the value at $s=0$ of the first derivative: the first one is in terms of $\theta_p$,
\begin{equation}\label{real}
L'\left(0,\sfrac{\cdot}{p}\right) =\log \theta_p,
\end{equation}
and the second one is in terms of the class number $h_p$ and the fundamental unit $\epsilon_p>1$,
\begin{equation}\label{CNF real}
L'\left(0,\sfrac{\cdot}{p}\right)=h_p\,\log \epsilon_p.
\end{equation}

Dirichlet also proved two formulas for an imaginary quadratic field. For simplicity, let $p\equiv 7\!\!\mod 8$ be a prime. Now the L-series $L\left(s,\sfrac{\cdot}{p}\right)$  does not vanish at $s=0$. His first formula states
\begin{equation}\label{im}
L\left(0,\sfrac{\cdot}{p}\right)=\sum_{0<a<p/2}\sfrac{a}{p},
\end{equation}
and the second one is in terms of the class number $h_{-p}$ of $\BQ[\sqrt{-p}]$
 \begin{equation}\label{CNF im}
L\left(0,\sfrac{\cdot}{p}\right)=h_{-p}.
\end{equation}
A non-trivial corollary is a finite expression for the class number 
\begin{equation*}
h_{-p}=\sum_{0<a<p/2}\sfrac{a}{p},
\end{equation*}
i.e., the difference between the number of square residues and of non-square residues in the interval $(0,p/2)$.

In 1952, Heegner discovered a way to construct rational points (over some number fields) on elliptic curves using special values of modular functions, in a manner similar to Dirichlet's solutions to Pell's equation. For instance, the elliptic curve
\begin{equation*}
E:y^2=x^3-1728
\end{equation*} 
is parameterized by modular functions $(\gamma_2,\gamma_3)$, where
\begin{equation*}
\gamma_2(z)=\frac{E_4}{\eta^8}=\frac{1+240\sum_{n=1}^\infty\sigma_{3}(n)q^n}{q^{1/3}\prod_{n= 1}^\infty (1-q^n)^{8}  }
,\,\quad
\gamma_3(z)=\frac{E_6}{\eta^{12}}=\frac{1-504\sum_{n=1}^\infty\sigma_{5}(n)q^n}{q^{1/2}\prod_{n= 1}^\infty (1-q^n)^{12}  }.
\end{equation*}
Here as customary, $z\in \CH$ is on the upper half plane, and $q=e^{2\pi i z}$. Then Heegner's strategy is to evaluate $(\gamma_2,\gamma_3)$ at $z$ in an auxiliary imaginary quadratic field $K$ to construct points of $E$ with coordinates in an abelian extension of $K$. The theorem of Gross and Zagier \cite{GZ} then relates the N\'eron--Tate heights of Heegner's points to special values of the {\em first order} derivatives of certain $L$-functions, providing an analog of Dirichlet's formula for a real quadratic field \eqref{real}. In this new context, an analog of Dirichlet's formula for an imaginary quadratic field \eqref{im} is the theorem of Waldspurger \cite{W85} relating toric period integrals (cf.\ \S\ref{sss:Wal}) to special values of L-functions of the same sort as in the work of Gross and Zagier.

 A natural question is to generalize the constructions of Dirichlet and of Heegner to higher dimensional algebraic varieties, and at the same time to generalize their relation to appropriate L-values. As a partial answer to this question, in this report we will consider some special algebraic cycles on Shimura varieties, cf.\ \S\ref{sss:pair Sh}. A class of such special cycles are the arithmetic diagonal cycles that appear in the arithmetic Gan--Gross--Prasad conjecture \cite[\S 27]{GGP}.

A parallel question is the generalization of the formulas of Dirichlet \eqref{im} and of Waldspurger to automorphic L-functions on a higher rank reductive group $\RG$ over a global field $F$. In this direction, we will consider automorphic period integrals for a spherical subgroup $\RH$ of $\RG$, i.e., 
$$
\int_{\RH(F)\bs\RH(\BA)}\phi(h)\,dh,
$$
where $\phi$ is a function on $\RG(F)\bs\RG(\BA)$, and $\BA$ the the ring of ad\`eles of $F$. The central object is the quotient $\RH(F)\bs\RH(\BA)$ sitting inside $\RG(F)\bs\RG(\BA)$. This paradigm may be viewed as a degenerate case of special cycles sitting inside the ambient Shimura variety. 

To study the relationship between automorphic periods and L-values, Jacquet invented the relative trace formula, cf.\ \cite{J05}.  In \cite{Z12} we adopted the relative trace formula approach to study height pairings of special cycles in the arithmetic Gan--Gross--Prasad conjecture. In this context, we formulated local conjectures (on intersection numbers of the arithmetic diagonal cycle on Rapoport--Zink spaces), namely the {\em arithmetic fundamental lemma} conjecture (cf.\ \cite{Z12}), and the {\em arithmetic transfer} conjecture by Rapoport, Smithling and the author (cf.\ \cite{RSZ1,RSZ2}).  In this report we will review the approach, the conjectures and the status.

Another natural question is in the direction of higher order derivatives of L-functions in the Gross--Zagier formula.   In  \cite{YZ1} Yun and the author found a geometric interpretation of the higher order derivative in the functional field setting, in terms of special cycles on the moduli stack of Drinfeld Shtukas with multiple number of modifications, cf.\ \S\ref{s:HGZ}. To generalize this to the number field case, one is led to the tantalizing question of finding Shtukas over number fields.

Limited by the length of this report, we will not discuss the analog of the class number formula \eqref{CNF real} and \eqref{CNF im} in our setting. In the context of elliptic curves over $\BQ$, this is the conjectural formula of Birch and Swinnerton--Dyer.  When the analytic rank is one (resp.\ zero), we have an equivalent statement in terms of the divisibility of Heegner points (resp.\ of normalized toric period integrals), thanks to the formula of Gross--Zagier (resp.\ of Waldspurger). Much of the equivalent statement has been proved in the past thirty years. Beyond elliptic curves (or modular forms on $\GL_2$), there have been many recent developments where special cycles play a crucial role in the study of Selmer groups of Bloch--Kato type. We hope to return to this topic on another occasion. 

\subsection*{Acknowledgments} The author thanks B.~Gross, M.~Rapoport, Y.~Sakellaridis, B.~Smithling and S.~Zhang for helpful comments.

\subsection*{Notation}
Let $F$ be a global field (unless otherwise stated), i.e., a number filed or a function field (of a geometrically connected smooth proper curve $X$ over a finite field $k=\BF_q$). Let $\BA=\BA_F=\prod'_{v}F_v$ be the ring of ad\`eles, the restricted direct product over all completions $F_v$  of $F$.  The ring of integers in a non-archimedean local field $F_v$ is denoted by $O_{F_v}$. For a subset of places $S$, we let $\BA_S=\prod_{v\in S}'F_v$. When $F$ is a number field, we write $\BA=\BA_f\times \BA_{\infty}$, where $\BA_f$ is the ring of the finite ad\`eles and $\BA_\infty=F\otimes_\BQ\BR$.

For a field extension $F'/F$ and an algebraic group $\RG$ over $F'$, we denote by $\RR_{F'/F}\RG$ the Weil restriction of scalars. We denote by $\BG_m$ the multiplicative group.

\section{Automorphic periods and L-values}

\subsection{Spherical pairs and automorphic periods}
\subsubsection{Automorphic quotient}\label{sss:aut quot}
Let $\RG$ be an algebraic group over $F$. We define the {\em automorphic quotient} associated to $\RG$ to be the quotient topological space
\begin{equation}
[\RG]\colon=\RG(F)\bs \RG(\BA).
\end{equation}
Let $K\subset \RG(\BA)$ be a subgroup.  When $F$ is a function field, we assume that $K$ is a compact open subgroup of $\RG(\BA)$. When $F$ is a number field, we assume that $K$ is a product $K_f\cdot K_\infty$ where $K_f\subset \RG(\BA_f)$  is a compact open subgroup, and $K_\infty$ is a suitable subgroup of $\RG(\BA_\infty)$.  We then define a quotient 
\begin{equation}
[\RG]_K\colon=\RG(F)\bs\left[ \RG(\BA)/K\right].
\end{equation}

When $F$ is a function field, this is a discrete set (or naturally as a groupoid). When $F=k(X)$, $\RG=\GL_n$, and $K=\prod_{v}\GL_n(O_{F_v})$, the groupoid $[\RG]_K$ is naturally isomorphic to the groupoid $\Bun_{n}(k)$, the $k$-points of $\Bun_{n}$  (the stack of vector bundles of rank $n$ on $X$).

\subsubsection{Automorphic period}\label{sss: per}
From now on let $\RG$ be a reductive group over $F$. Let $\RH\subset \RG$ be a subgroup. Let $Z_\RG$ be the center of $\RG$ and let $Z=\RH\cap Z_\RG$. Let $\CA_0(\RG)$ be the space of cuspidal automorphic forms on $[\RG]$, invariant under the action of $Z(\BA)$. Then the {\em automorphic $H$-period integral} is defined by
\begin{align*}
  \xymatrix@R=0ex{
	  \sP_{\RH} \colon \CA_0(\RG)\ar[r]  &  \BC  \\
		\phi \ar@{|->}[r]
		   &  \int_{Z(\BA)\bs [\RH]}\phi(h)\,dh.
	}
	\end{align*}
\begin{remark}The name ``automorphic period" is different from the period in the context of comparison theorems between various cohomology theories. However, some special cases of the automorphic period integrals may yield periods in the de Rham--Betti comparison theorem. For instance, this happens when the integral can be turned into the form $\int_{Z(\BA)\bs [\RH]_{K_\RH}}\omega_\phi $ for a closed differential form $\omega_{\phi}$ on the real manifold $[\RG]_K$ for suitable $K$, and $K_\RH=K\cap \RH(\BA)$. 
\end{remark}

Let
$\pi$ be a (unitary throughout the article) cuspidal automorphic representation of $\RG(\BA)$. For simplicity, we assume that there is a unique embedding $\pi\incl \CA_0(\RG)$  (in many applications this is the case). We consider the restriction of $\sP_\RH$ to $\pi$; this defines an element  in $\sP_{\RH,\pi}\in\Hom _{\RH(\BA)}(\pi, \BC)$. 
\begin{defn}[Jacquet]
A cuspidal automorphic representation  $\pi$ is {\it (globally) distinguished by $\RH$} if the linear functional $\sP_{\RH,\pi}\in \Hom _{\RH(\BA)}(\pi, \BC)$ does not vanish, i.e., there exists some $\phi\in \pi$ such that  $\sP_{\RH}(\phi)\neq 0$.
\end{defn}

It is also natural to consider a twisted version. Let $\chi$ be a character of $Z(\BA)\RH(F)\bs \RH(\BA)$. We define the automorphic $(\RH,\chi)$-period integral $\sP_{\RH,\chi}$ in a similar manner,
$$
\sP_{\RH,\chi}(\phi)=\int_{Z(\BA)\bs [\RH]}\phi(h)\chi(h)\,dh.
$$

If $\pi$ is distinguished by $\RH$, then $\Hom _{\RH(\BA)}(\pi, \BC)\neq 0$, and in particular,  $\Hom _{\RH(F_v)}(\pi_v, \BC)\neq 0$ for every place $v$. We then say that $\pi_v$ is (locally) distinguished by $\RH(F_v)$ if $\Hom _{\RH(F_v)}(\pi_v, \BC)\neq 0$. Very often the automorphic period integral $\sP_\RH$ behaves nicely only when the pair $(\RH,\RG)$ satisfies certain nice properties, such as 
\begin{altenumerate}
\item the multiplicity-one property  $\dim \Hom _{\RH(F_v)}(\pi_v, \BC) \leq 1$ holds for all $v$, or the multiplicity can be described in a certain nice way, and
\item the locally distinguished representations can be characterized in terms of L-parameters.  
\end{altenumerate}

Jacquet and his school have studied the distinction for many instances (locally and globally, cf.\ \cite{J05}). A large class of $(\RH,\RG)$ called {\em spherical pairs} are expected to have the above nice properties, according to the work of Sakellaridis \cite{Sa08} and his joint work with Venkatesh \cite{SV}. 

\subsubsection{Spherical pairs}
\label{sss:sph}
Let $(\RH,\RG)$ be over a field $F$ (now arbitrary). The reductive group $\RG$ acts on $\RG/\RH$ by left multiplication. If $F$ is an algebraically closed field, we say that the pair $(\RH,\RG)$ is {\em spherical} if a Borel subgroup $B$ of $\RG$ has an open dense orbit on $\RG/\RH$ \cite{Sa08}. Over a general field $F$, the pair $(\RH,\RG)$ is said to be {\em spherical} if its base change to an algebraic closure of $F$ is spherical. We then call $\RH$ a  spherical subgroup of $\RG$.

Here are some examples.
\begin{altenumerate}\label{eg sph}
\item Whittaker pair $(\RN,\RG)$, where $\RG$ is quasi-split, and $\RN$ is a maximal nilpotent subgroup. 
 \item The pair $(\RG,\RG\times\RG)$ with the diagonal embedding.
\item Symmetric pair $(\RH,\RG)$, where $\RH$ is the fixed point locus of an involution $\theta:\RG\to \RG$. This constitutes a large class of spherical pairs, including 
 \begin{altitemize}
 \item the unitary periods of Jacquet and Ye (cf.\ the survey \cite{Of09}), where $\RG=\RR_{F'/F}\GL_n$, and $\RH$ is a unitary group attached to a Hermitian space with respect to a quadratic extension $F'/F$. This period is related to the quadratic base change for the general linear group.
 \item the Flicker--Rallis periods \cite[\S3.2]{Z14b},  where $\RG=\RR_{F'/F}\GL_n$, and $\RH=\GL_{n}$, for quadratic $F'/F$. This period is related to the quadratic base change for the unitary group.

 \item the linear periods of Friedberg--Jacquet \cite{FJ}, where  $\RH=\GL_{n/2}\times \GL_{n/2}$ embeds in  $\RG= \GL_{n}$ by $$(a,d)\mapsto\diag(a,d).$$
 \end{altitemize}
\item the Rankin--Selberg pair (named after its connection to Rankin--Selberg convolution L-functions), where $\RG=\GL_{n-1}\times \GL_{n}$, and $\RH=\GL_{n-1}\incl G$ with the embedding
 $$g\mapsto (g,\,\diag(g,1)).$$ (There are also spherical subgroups of $\GL_{m}\times \GL_n$ when $|n-m|>1$ involving non-reductive subgroups.) 
\item the Gan--Gross--Prasad pairs $(\SO_{n-1},\SO_{n-1}\times\SO_n)$ and $(\RU_{n-1},\RU_{n-1}\times\RU_n)$; see \S\ref{sss:GGP}. They resemble the Rankin--Selberg pairs. (There are also spherical subgroups of $\SO_{m}\times\SO_n$ and $\RU_{m}\times\RU_n$ when $|n-m|>1$ involving non-reductive subgroups).
\end{altenumerate}

\subsection{The global Gan--Gross--Prasad conjecture and the Ichino--Ikeda refinement}
 
Let $F$ be a number field for the rest of this section.  
\subsubsection{Waldspurger formula}\label{sss:Wal}
 Let $B$ be a quaternion algebra over $F$ and let $\RG=B^\times$ (as an $F$-algebraic group).  Let $F'/F$ be a quadratic extension of number fields and denote by $\RT$ the torus ${\rm R}_{F'/F}\BG_m$. Let $F'\incl B$ be an embedding of $F$-algebras, and $\RT\incl \RG$ the induced embedding of $F$-algebraic groups. Then $(\RT,\RG)$ is a spherical pair. In \cite{W85}, Waldspurger studied the automorphic period integral $\sP_{T,\chi}$ (sometimes called the {\em toric} period), and he proved an exact formula relating the square $|\sP_{T,\chi}|^2$ to a certain central L-value.

Below we will consider one generalization of Waldspurger's formula to higher rank groups, i.e., the global conjectures of Gan--Gross--Prasad and Ichino--Ikeda.  For this report, we implicitly assume the endoscopic classification of Arthur for orthogonal and unitary groups.

\subsubsection{The global Gan--Gross--Prasad conjecture}\label{sss:GGP}
 
Gan, Gross, and Prasad \cite{GGP} proposed a series of precise conjectures regarding the local and global distinction for $(\RH,\RG)$ when $\RG$ is a classical group (orthogonal, unitary and symplectic),  extending the conjectures of Gross and Prasad \cite{GP1,GP2} for orthogonal groups. 
 
We recall their global conjectures in orthogonal and Hermitian cases. For simplicity, we restrict to the case when the spherical subgroup is reductive. Let $F$ be a number field, and let $F'=F$ in the orthogonal case and $F'$ a quadratic extension of $F$ in the Hermitian case. Let $W_n$ be a non-degenerate orthogonal space or Hermitian space with $F'$-dimension $n$.  Let $W_{n-1}\subset W_n$ be a non-degenerate subspace of codimension one. Let $\RG_{i}$ be $\SO(W_i)$ or $\RU(W_i)$ for $i=n-1,n$, and $\delta:\RG_{n-1}\incl\RG_n$ the induced embedding. Let 
\begin{align}
 \RG=\RG_{n-1}\times \RG_{n},
	\quad\,\quad
	\RH=\RG_{n-1},
\end{align} with the ``diagonal" embedding $\Delta:\RH \incl \RG$ (i.e., the graph of $\delta$). The pair $(\RH,\RG)$ is spherical and we call it the {\em Gan--Gross--Prasad pair}.

Let $\pi=\pi_{n-1}\boxtimes \pi_{n}$ be a tempered cuspidal automorphic
representation of $\RG(\BA)$. The central L-values of certain automorphic L-functions $L(s,\pi,R)$ show up in their conjecture, where $R$ is a finite dimensional representation of the L-group $^L\RG$, cf.\ \cite[\S22]{GGP}. We can describe the L-function as the Rankin--Selberg convolution of suitable automorphic representations on general linear groups. For $i\in\{n-1,n\}$, let $\Pi_{i,F'}$ be the endoscopic functoriality transfer of $\pi_i$ from $\RG_i$ to suitable $\GL_N(\BA_{F'})$: in the Hermitian case, this is the base change of $\pi_i$ to
$\GL_i(\BA_{F'})$; and in the orthogonal case, this is the endoscopic transfer from $\RG_i(\BA)$ to $\GL_i(\BA)$ (resp.\ $\GL_{i-1}(\BA)$) if $i$ is even (resp.\ odd). Then the L-function $L(s,\pi,R)$ can be defined more explicitly as the Rankin--Selberg convolution L-function $L(s,\Pi_{n-1,F'}\times
\Pi_{n,F'})$.

We are ready to state the global Gan--Gross--Prasad conjecture  \cite[\S24]{GGP}.
\begin{conj}\label{conj GGP}
Let $\pi$ be a tempered cuspidal automorphic
representation of $\RG(\BA)$.
 The following statements are equivalent.
\begin{altenumerate}
\item The automorphic $\RH$-period integral does not vanish on $\pi$, i.e., $\sP_{\RH}(\phi)\neq 0$ for some $\phi\in\pi$.
\item The space $\Hom_{\RH(\BA)}(\pi,\BC)\neq 0$ and the central value $L(\frac{1}{2},\pi,R)\neq 0$.
\end{altenumerate}
\end{conj}

\begin{remark} It is known that the pair $(\RH(F_v),\RG(F_v))$ satisfies the multiplicity one property by Aizenbud, Gourevitch, Rallis, and Schiffmann \cite{AGRS} for $p$-adic local fields, and by Sun and Zhu \cite{SZ} for archimedean local fields. The local conjectures of Gan, Gross, and Prasad \cite[\S17]{GGP} specify the member $\pi_v$ in a generic Vogan L-packet (cf.\ \cite[\S9-11]{GGP}) with $\dim\Hom_{\RH(F_v)}(\pi_v,\BC)=1$, in terms of local root numbers associated to the L-parameter.  Their local conjectures are mostly proved by M\oe glin and Waldspurger \cite{W12,MW} ($p$-adic orthogonal groups), and Beuzart-Plessis \cite{BP15a,BP15b} (unitary groups over $p$-adic and archimedean local fields).
\end{remark}

We have the following result.
\begin{thm}
Let $\RG=\RU(W_{n-1})\times\RU(W_n)$ for Hermitian spaces $W_{n-1}\subset W_n$ over a quadratic extension $F'$ of $F$.
 Let $\pi$ be a tempered cuspidal automorphic representation of $\RG(\BA)$ such that  there exist a non-archimedean place $v$ of $F$ split in $F'$ where $\pi_{v}$ is supercuspidal. Then Conjecture \ref{conj GGP} holds.
\end{thm}

This was proved in \cite{Z14a} under a further condition on the archimedean places, which was later removed by Xue \cite{Xue}. The local condition above was due to a simple version of the Jacquet--Rallis relative trace formula in \cite{Z14a}, cf.\ \S\ref{ss:JR}. Zydor and Chaudouard \cite{CZ,Zy} have made progress towards the full relative trace formula that should remove the local condition.

\begin{remark}  Ginzburg, Jiang, and Rallis \cite{GJR1,GJR2} have proved the direction ${\rm (i)\imp (ii)}$ of Conjecture \ref{conj GGP} for both the orthogonal and Hermitian cases when the group $\RG$ is quasi-split  and the representation $\pi$ is (globally) generic.
\end{remark}

\subsubsection{Ichino--Ikeda refinement}
For many applications, we would like to have a refined version of the Gan--Gross--Prasad conjecture, analogous to the Waldspurger formula for the toric period $\sP_{\RT,\chi}$ in \S\ref{sss:Wal}. We recall the refinement of Ichino and Ikeda \cite{II} (for orthogonal groups; later their idea was carried out for unitary groups by N. Harris in \cite{HN}). 

Let $L(s,\pi, \Ad)$  be the adjoint L-function (cf.\ \cite[\S7]{GGP}).
Denote $\Delta_{n}=L(M^{\vee}(1))$ where $M^\vee$ is the motive
dual to the motive $M$ associated to $\RG_{n}$ by
Gross \cite{G97}. It is a product of special values of Artin L-functions.
We will be interested in the following combination of L-functions,
\begin{align}
\sL(s,\pi)=\Delta_{n}\frac{L(s,\pi,R)}{L(1,\pi, \Ad)}.
\end{align}
We also write $\sL(s,\pi_v)$ for the corresponding local factor at a place $v$.

Let $\pi_v$ be an irreducible tempered unitary representation of $\RG(F_v)$ with an invariant inner product $\langle \cdot,\cdot\rangle _v$.
Ichino and Ikeda construct a canonical element in the space $\Hom_{\RH(F_v)}(\pi_v,\BC)\otimes\Hom_{\RH(F_v)}(\ov\pi_v,\BC)$ by integrating matrix coefficients: for $\phi_v,\varphi_v\in \pi_v$,  
\begin{align}\label{alpha}
\wt\alpha_v(\phi_v,\varphi_v)=\int_{\RH(F_v)}{\langle \pi_v(h)\phi_v,\varphi_v
\rangle _v}\, dh.
\end{align}Ichino and Ikeda  showed that the integral converges absolutely for all (tempered) $\pi$. In fact, the convergence holds for any {\em strongly tempered} pair $(\RH,\RG)$, cf.\ \cite{SV}. When $\pi_v$ is unramified and the vectors $\phi_v,\varphi_v$ are fixed by a hyperspecial compact open $\RG(O_{F_v})$ such that $\langle \phi_v,\varphi_v \rangle_v=1$, we have 
$$
\wt\alpha_v(\phi_v,\varphi_v)=\sL(\frac{1}{2},\pi_v) \cdot \vol(\RH(O_{F_v})).$$
We normalize the local canonical invariant form:
\begin{align}\label{alpha1}
\alpha_v(\phi_v,\varphi_v)=\frac{1}{\sL(\frac{1}{2},\pi_v)}\wt\alpha_v(\phi_v,\varphi_v).
\end{align}

We endow $\RH(\BA)$  (resp.\ $\RG(\BA)$) with their Tamagawa measures and
$[\RH]$ (resp.\ $[\RG]$) with the quotient measure by the counting measure on $\RH(F)$ (resp.\ $\RG(F)$).
We choose the Haar measure $dh$ on $\RH(\BA)$ and the measures $dh_v$ on $\RH(F_v)$ such that
$dh=\prod_{v}dh_v.$ Let $\pair{\phi,\varphi}_{\rm Pet}$ be the Petersson inner product of $\phi,\varphi\in\pi=\otimes_v\pi_v$, and choose the local inner products $\langle \cdot,\cdot\rangle _v$ such that $\pair{\phi,\varphi}=\prod_v \pair{\phi_v,\varphi_v}_v$ for $\phi=\otimes_v\phi_v$ and $\varphi=\otimes_v\varphi_v$. 

We can now state the Ichino--Ikeda conjecture \cite[Conj.\ 1.5 and 2.1]{II} that refines the global Gan--Gross--Prasad conjecture \ref{conj GGP}.

\begin{conj}\label{conj IIH}
Let $\pi$ be a tempered cuspidal automorphic
representation of $\RG(\BA)$.
Then for $\phi=\otimes_{v}
\phi_v\in \pi$,
\begin{equation}\label{fml *1}
\big|\sP_\RH(\phi)\big|^2=2^{-\beta_\pi}\sL(1/2,\pi)
\prod_{v}\alpha_v(\phi_v,\phi_v),
\end{equation}
where $\beta_\pi$ is the rank of a finite elementary $2$-group associated to the L-parameter of $\pi$.
\end{conj}

\begin{remark}
If $\Hom_{\RH(\BA)}(\pi,\BC)=0$, both sides of \eqref{fml *1} vanish. 
 \end{remark}

\begin{remark}In the orthogonal case and $n=3$, the refined conjecture is exactly the same as the formula of Waldspurger. When $n=4$ the conjecture is proved by Ichino \cite{Ich}. Little is known in the higher rank case, cf.\ the survey \cite{Gan13}. \end{remark}

In the Hermitian case and $n=2$, the conjecture follows from Waldspurger's formula \cite{HN}. In general, 
we have the following result, due to the author \cite{Z14b} and Beuzart-Plessis \cite{BP16}. \begin{thm} Let $\RG=\RU(W_{n-1})\times\RU(W_n)$ for Hermitian spaces $W_{n-1}\subset W_n$ over a quadratic extension $F'$ of $F$.  Let $\pi$ be a tempered cuspidal automorphic representation of $\RG(\BA)$. Assume that 
\begin{altenumerate}
\item
there exists a non-archimedean place $v$ of $F$ split in $F'$ such that $\pi_{v}$ is supercuspidal, and 
\item all archimedean places of $F$ are split in $F'$.
\end{altenumerate}
 Then Conjecture \ref{conj IIH} holds.
\end{thm}For some ingredients of the proof, see \S\ref{sss:JR sph}.

\subsubsection{Reformulation in terms of spherical characters}
\label{sss:sph GGP}

Let $\pi$ be a cuspidal automorphic representation of $\RG(\BA)$. We define the {\em global spherical character} $\BI_\pi$ as the distribution on $\RG(\BA)$
\begin{align}\label{global char}
\BI_\pi(f):=\sum_{\phi\in {\rm OB}(\pi)} \sP_\RH(\pi(f)\phi)\ov{\sP_\RH(\phi)},\quad f\in \sC^\infty_c(\RG(\BA)),
\end{align}
where the sum runs over an orthonormal basis $ {\rm OB}(\pi)$ of $\pi$ (with respect to the Petersson inner product). Note that $\BI_\pi$ is an {\em eigen-distribution} for the spherical Hecke algebra $\sH^S(\RG)$ away from a sufficiently large set $S$ (including all bad primes), in the sense that for all $f=f_S\otimes f^S$ with $f^S\in \sH^S(\RG)$ and $f_S\in  \sC_c^\infty(\RG(\BA_S))$,
\begin{align}\label{eigen dis}
\BI_\pi(f)=\lambda_{\pi^S}\left(f^S\right)\,\BI_\pi\left(f_S\otimes 1_{K^S}\right),
\end{align}
where $\lambda_{\pi^S}$ is the ``eigen-character" of $\sH^S(\RG)$ associated to $\pi^S$.

 We define  the {\em local spherical character} in terms of the local canonical invariant form $\alpha_v$ in \eqref{alpha1},
\begin{align}\label{def local char G}
\BI_{\pi_v}(f_v):=\sum_{\phi_v\in {\rm OB}(\pi_v)} \alpha_v(\pi_v(f_v)\phi_v,\phi_v),\quad f_v\in \sC^\infty_c(\RG(F_v)),
\end{align}
where the sum runs over an orthonormal basis $ {\rm OB}(\pi_v)$ of $\pi_v$.

In \cite[Conj. 1.6]{Z14b} the author stated an alternative version of the Ichino--Ikeda conjecture in terms of spherical characters.
\begin{conj}\label{conj dis}
Let $\pi$ be a tempered cuspidal automorphic representation of $\RG(\BA)$.
Then for all pure tensors $f=\bigotimes_{v} f_v\in \sC^\infty_c(\RG(\BA))$, 
$$
\BI_\pi(f)=2^{-\beta_\pi}\sL(1/2,\pi) \prod_{v}\BI_{\pi_v}(f_v).
$$
\end{conj} 
By \cite[Lemma 1.7]{Z14b}, this conjecture is equivalent to Conjecture \ref{conj IIH}. This new formulation is more suitable for the relative trace formula approach, cf.\ \S\ref{sss:JR sph}. This also inspires us to state a version of the refined arithmetic Gan--Gross--Prasad conjecture where the Ichino--Ikeda formulation does not seem to apply directly, cf.\ \S\ref{sss: sph AGGP}, Conjecture \ref{conj:sph AGGP}.

\section{Special cycles and L-derivatives}

\subsection{Special pairs of Shimura data and special cycles}
\subsubsection{}\label{sss:pair Sh}
To describe our set-up, we introduce the concept of {\em a special pair of Shimura data}. Let  $\BS$ be the torus $\RR_{\BC/\BR}\BG_m$ over $\BR$ (i.e., we view $\BC^\times$ as an $\BR$-group). Recall that a Shimura datum $\bigl(\RG,X_\RG\bigr)$ consists of a reductive group $\RG$ over $\BQ$, and a $\RG(\BR)$-conjugacy class $X_\RG=\{h_\RG\}$ of $\BR$-group homomorphisms $h_\RG:\BS\to \RG_\BR$ (sometimes called Shimura homomorphisms)  satisfying Deligne's list of axioms \cite[1.5]{De}. In particular, $X_\RG$ is a Hermitian symmetric domain. 
 \begin{defn}
 A {\it special pair} of Shimura data is a homomorphism  \cite[1.14]{De} between two Shimura data 
 $$
 \xymatrix{\delta\colon\bigl(\RH,X_\RH\bigr)\ar[r]& \bigl(\RG,X_\RG\bigr) }
 $$ such that
  \begin{altenumerate}
 \item the homomorphism $\delta: \RH\to \RG$ is injective such that the pair $(\RH,\RG)$ is spherical, and 
 \item the dimensions of $X_\RH$ and $X_\RG$ (as complex manifolds) satisfy 
 $$
\dim_\BC X_\RH=\bigg\lfloor \frac{\dim_\BC X_\RG}{2}\bigg\rfloor.
 $$ 
 \end{altenumerate}
 \end{defn} 
 
In particular, we enhance a spherical pair $(\RH, \RG)$ to a homomorphism of Shimura data $(\RH,X_\RH)\to (\RG,X_\RG)$. 

 \begin{remark}
 It seems an interesting question to enumerate special pairs of Shimura data.  In fact, we may consider the analog of special pairs of Shimura data in the context of {\em local Shimura data} \cite{RV}. It seems more realistic to enumerate the pairs in the local situation.
\end{remark}

For a Shimura datum $  (\RG,X_\RG)$ we have a projective system of {\em Shimura varieties} $\{\Sh_{K}(\RG)\}$, indexed by compact open subgroups $K\subset\RG(\BA_f)$,  of smooth quasi--projective varieties (for neat $K$) defined over a number field $E$---the reflex field of $(\RG,X_\RG)$.

For a special pair of Shimura data $(\RH,X_\RH)\to (\RG,X_\RG)$, compact open subgroups $K_\RH\subset \RH(\BA_f)$ and $K_\RG\subset\RG(\BA_f)$ such that $K_\RH\subset K_\RG$, we have a finite morphism (over the reflex field of $(\RH,X_\RH)$)
$$
\delta_{K_\RH,K_\RG}:\Sh_{K_\RH}(\RH)\to \Sh_{K_\RG}(\RG).
$$
The cycle $z_{K_\RH,K_\RG}:=\delta_{K_\RH,K_\RG,\ast} [\Sh_{K_\RH}(\RH)]$ on $ \Sh_{K_\RG}(\RG)$ will be called the {\em special cycle} (for the level $(K_\RH,K_\RG)$).  Very often we choose $K_\RH=K_\RG\cap \RH(\BA_f)$ in which case we simply denote the special cycle by $z_{K_\RG}$.
 \begin{remark}
 \begin{altenumerate}
 \item
Note that here our special cycles are different from those appearing in \cite{KRY,KR-U2}.
\item When $\dim_\BC X_\RG$ is even, the special cycles are in the middle dimension. When $\dim_\BC X_\RG$ is odd, the special cycles are just below the middle dimension, and we will say that they are in the {\em arithmetic  middle dimension} (in the sense that, once extending both Shimura varieties to suitable integral models, we obtain cycles in the middle dimension). 
\end{altenumerate}
 \end{remark}

The special cases in the middle dimension are very often related to the study of Tate cycles and automorphic period integrals, e.g., in the pioneering example of Harder, Langlands, and Rapoport \cite{HLR}, and many of its generalizations.

Below we focus on the case where the special cycles are in the arithmetic  middle dimension. 
 \subsubsection{Gross--Zagier pair}
  In the case of the Gross--Zagier formula \cite{GZ}, one considers an embedding of an {\em imaginary} quadratic field $F'$ into $\Mat_{2,\BQ}$ (the algebra  of $2\times 2$-matrices), and the induced embedding  
  $$\RH=\RR_{F'/\BQ}\BG_m\incl \RG=\GL_{2,\BQ}.$$   Note that $\RH_\BR\simeq\BC^\times$ as $\BR$-groups (upon a choice of embedding $F'\incl \BC$). This defines  $h_\RH: \BS\to \RH_\BR$, and its composition with the embedding $ \RH_\BR\to \RG_\BR$ defines $h_{\RG}:\BS\to \BG_\BR$. We obtain a  special pair $(\RH,X_\RH)\to (\RG,X_\RG)$, 
where
  $$
  \dim  X_\RG=1,\quad\quad \dim  X_\RH=0.
  $$
  
  In the general case, we replace $F'/\BQ$ by a CM extension $F'/F$ of a totally real number field $F$, and replace $\Mat_{2,\BQ}$ by a quaternion algebra $B$ over $F$ that is ramified at all but one archimedean places of $F$. X. Yuan, S. Zhang, and the author proved Gross--Zagier formula  in this generality in \cite{YZZ1}.

 \subsubsection{Gan--Gross--Prasad pair}
  \label{sss:GGP pair}
 Let $(\RH,\RG)$ be the Gan--Gross--Prasad pair in \S\ref{sss:GGP}, but viewed as algebraic groups over $\BQ$ (i.e., the pair $(\RR_{F/\BQ}\RH,\RR_{F/\BQ}\RG)$). The groups are associated to an embedding  $W_{n-1}\subset W_{n}$ of orthogonal or Hermitian spaces (with respect to $F'/F$). Now we impose the following conditions.
 \begin{altenumerate}
 \item $F$ is a totally real number field, and in the Hermitian case $F'/F$ is a CM (=totally imaginary quadratic) extension.
 \item  For an archimedean place $\varphi\in\Hom(F,\BR)$, denote by $\sgn_{\varphi}(W)$ the signature of $W\otimes_{F,\varphi}\BR$ as an orthogonal or Hermitian space over $F'\otimes_{F,\varphi}\BR$. Then there exists a  distinguished real place $\varphi_0\in\Hom(F,\BR)$ such that
 $$
 \sgn_{\varphi}(W_{n})=\begin{cases}(2,n-2),&\varphi=\varphi_0
 \\(0,n),& \varphi\in\Hom(F,\BR)\setminus\{\varphi_0\}
 \end{cases}
 $$
 in the orthogonal case, and 
  $$
 \sgn_{\varphi}(W_{n})=\begin{cases}(1,n-1),&\varphi=\varphi_0
 \\(0,n),& \varphi\in\Hom(F,\BR)\setminus\{\varphi_0\}
 \end{cases}
 $$in the Hermitian case. In addition, the quotient $W_n/W_{n-1}$ is negative definite at every $\varphi\in\Hom(F,\BR)$ (so the signature of $W_{n-1}$ is given by similar formulas).
\end{altenumerate} Then Gan, Gross, and Prasad \cite[\S27]{GGP} prescribe Shimura data that enhance
the embedding $\RH\incl \RG$ to a homomorphism of Shimura data $(\RH,X_\RH)\to (\RG,X_\RG)$, where the dimensions are 
$$\begin{cases}
\dim  X_\RG=2n-5,\quad \dim  X_\RH=n-3,\quad \text{in the orthogonal case,}\\
\dim  X_\RG=2n-3,\quad \dim  X_\RH=n-2,\quad \text{in the Hermitian case}.
\end{cases}
$$

\subsection{The arithmetic Gan--Gross--Prasad conjecture}

\subsubsection{Height pairings}
\label{sss std conj}

 Let $X$ be a smooth proper variety over a number field $E$, and let $\Ch^i(X)$ be the group of codimension-$i$ algebraic cycles on $X$ modulo rational equivalence. We have a cycle class map
$$
   \cl_i\colon \Ch^i(X)_\BQ\to H^{2i}(X),
$$where $H^{2i}(X)$ is the Betti cohomology
$   H^*\bigl(X(\BC),\BC\bigr).$
The kernel is the group of cohomologically trivial cycles, denoted by $\Ch^i(X)_{0}$.

 Conditional on some standard conjectures on algebraic cycles, there is a height pairing defined by Beilinson and Bloch, \begin{equation}
\label{eqn BB}
  \sform_\mathrm{BB}\colon \Ch^i(X)_{\BQ,0}\times \Ch^{d+1-i}(X)_{\BQ,0}\to \BR,\quad d=\dim X.
\end{equation}
This is unconditionally defined when $i=1$ (the N\'eron--Tate height), or when $X$ is an abelian variety \cite{Ku-ens}. 
In some situations, cf.\ \cite[\S6.1]{RSZ3}, one can define the height pairing unconditionally in terms of the arithmetic intersection theory of 
Arakelov and Gillet--Soul\'e \cite[\S4.2.10]{GS}. This is the case when
there exists a smooth proper model $\CX$ of $X$ over $O_E$ (this is also true for Deligne--Mumford (DM) stacks $X$ and $\CX$).

\subsubsection{The arithmetic Gan--Gross--Prasad conjecture}

We consider the special cycle in the Gan--Gross--Prasad setting \S\ref{sss:GGP pair}, which we also call the {\em arithmetic diagonal cycle} \cite{RSZ3}. We will state a version of the  arithmetic Gan--Gross--Prasad conjecture assuming some standard conjectures on algebraic cycles (cf.\ \cite[\S6]{RSZ3}), in particular, that we have the height pairing \eqref{eqn BB}. 

For each $K\subset \RG(\BA_f)$, one can construct ``Hecke--Kunneth" projectors that project the total cohomology of the Shimura variety $\Sh_{K}(\RG)$ (or its toroidal compactification) to the odd-degree part (cf.\ \cite[\S6.2]{RSZ3} in the Hermitian case; the same proof works in the orthogonal case).  Then we apply this projector to define a cohomologically trivial cycle $z_{K,0}\in \Ch^{n-1}\bigl(\Sh_{K}(\RG)\bigr)_{0}$ (with $\BC$-coefficient). The classes $\{z_{K,0}\}_{K\subset \RG(\BA_f)}$ are independent of the choice of our projectors (cf.\ \cite[Remark 6.11]{RSZ3}), and they form a projective system (with respect to push-forward).  

We form the colimit $$\Ch^{n-1}\bigl(\Sh(\RG)\bigr)_{0}:=\varinjlim_{K\subset \RG(\BA_f)} \Ch^{n-1}\bigl(\Sh_{K}(\RG)\bigr)_{0}.$$ The height pairing with $\{z_{K,0}\}_{K\subset \RG(\BA_f)}$ defines a linear functional
$$
\xymatrix{\sP_{\Sh(\RH)}\colon \Ch^{n-1}\bigl(\Sh(\RG)\bigr)_{0}\ar[r]&\BC}.
$$
  This is the arithmetic version of the automorphic period integral in \S\ref{sss: per}.   The group $ \RG(\BA_f)$ acts on the space $ \Ch^{n-1}\bigl(\Sh(\RG)\bigr)_{0}$. For any representation $\pi_f$ of $ \RG(\BA_f)$, let $\Ch^{n-1}\bigl(\Sh(\RG)\bigr)_{0}[\pi_f]$ denote the $\pi_f$-isotypic component of the Chow group $\Ch^{n-1}\bigl(\Sh(\RG)\bigr)_{0}$.
  
  We are ready to state the arithmetic Gan--Gross--Prasad conjecture \cite[\S27]{GGP}, parallel to Conjecture \ref{conj GGP}.
\begin{conj}\label{conj AGGP}
Let $\pi$ be a tempered cuspidal automorphic
representation of $\RG(\BA)$, appearing in the cohomology $H^{\ast}(\Sh(\RG))$.
 The following statements are equivalent.
\begin{altenumerate}
\item The linear functional $\sP_{\Sh(\RH)}$ does not vanish on the $\pi_f$-isotypic component $\Ch^{n-1}\bigl(\Sh(\RG)\bigr)_{0}[\pi_f]$.
\item  The space $\Hom_{\RH(\BA_f)}(\pi_f,\BC)\neq 0$ and the first order derivative  $L'(\frac{1}{2},\pi,R)\neq 0$.
\end{altenumerate}
\end{conj}

 In the orthogonal case with $n\leq 4$, and when the ambient Shimura variety is a curve ($n=3$), or a product  of three curves ($n=4$),  the conjecture is unconditionally formulated.  The case $n=3$ is proved by X. Yuan, S. Zhang, and the author in \cite{YZZ1}; in fact we proved a refined version.  When $n=4$ and in the triple product case (i.e., the Shimura variety $\Sh_{K}(\RG)$ is a product of three curves), X. Yuan, S. Zhang, and the author formulated a refined version of the above conjecture and proved it in some special cases, cf.\ \cite{YZZ2}.

\subsubsection{Reformulation in terms of  spherical characters}
\label{sss: sph AGGP}

In the Hermitian case, Rapoport, Smithling, and the author in \cite{RSZ3} stated a version of the arithmetic Gan--Gross--Prasad conjecture that does not depend on standard conjectures on algebraic cycles.

In fact we work with a variant of the Shimura data defined by Gan, Gross, and Prasad \cite[\S27]{GGP}. We modify the groups $\RR_{F/\BQ}\RG$ and $\RR_{F/\BQ}\RH$ defined previously
{\allowdisplaybreaks
\begin{align*}\RZ^\BQ&:=\GU_1=\bigl\{\,   z\in  \RR_{F'/\BQ}\BG_m \bigm| \Nm_{F'/F}(z) \in \BG_m  \,\bigr\} ,\\
\wt \RH &:={\rm G\left(U_1\times U(W_{n-1})\right)}=\bigl\{\, (z,h) \in  \RZ^\BQ\times \GU(W_{n-1})\bigm| \Nm_{F'/F}(z) = c(h) \,\bigr\},\\
\wt \RG &:={\rm G\left(U_1\times U(W_{n-1})\times U(W_{n}) \right)}\\
&\phantom{:}= \bigl\{\, (z,h,g) \in \RZ^\BQ \times \GU(W_{n-1})\times  \GU(W_{n})\bigm| \Nm_{F'/F}(z) = c(h)=c(g) \,\bigr\},
\end{align*}
where the symbol $c$ denotes the unitary similitude factor. 
Then we have 
\begin{equation}\label{proddec}
   \begin{gathered}
   \begin{gathered}
	\xymatrix@R=0ex{
      \wt \RH \ar[r]^-\sim  &  \RZ^\BQ \times \RR_{F/\BQ} \RH
	}
   \end{gathered},
	\qquad
   \begin{gathered}
	\xymatrix@R=0ex{
	   \wt \RG \ar[r]^-\sim  & \RZ^\BQ \times \RR_{F/\BQ} \RG	}
	\end{gathered}.
   \end{gathered}
\end{equation}

We then define natural Shimura data $\big(\wt \RH,\{h_{\wt \RH}\} \big)$ and $\big(\wt \RG,\{h_{\wt \RG}\}\big)$, cf.\ \cite[\S3]{RSZ3}. This variant has the nice feature that the Shimira varieties are of PEL type, i.e., the canonical models are related to moduli problems of abelian varieties with polarizations, endomorphisms, and level structures, cf.\ \cite[\S4--\S5]{RSZ3}. 

For suitable Hermitian spaces and a special level structure $ K_{\wt \RG}^\circ\subset \wt\RG(\BA_f)$, we can even define {\em smooth} integral models (over the ring of integers of the reflex field) of  the Shimura variety $\Sh_{ K_{\wt \RG}^\circ}(\wt \RG)$. For a general CM extension $F'/F$, it is rather involved to state this level structure \cite[Remark 6.19]{RSZ3} and define the integral models \cite[\S5]{RSZ3}. For simplicity, from now on we consider a special case, when $F=\BQ$ and $F'=F[\varpi]$ is an imaginary quadratic field. We further assume that the prime $2$ is split in $F'$. We choose $\varpi\in F'$ such that $(\varpi)\subset O_{F'}$ is the product of all ramified prime ideals in $O_{F'}$.  

We first define an auxiliary moduli functor $\CM_{(r,s)}$ over $\Spec O_{F'}$ for $r+s=n$ (similar to \cite[\S13.1]{KR-U2}). 
For a locally noetherian scheme $S$ over $\Spec O_{F'}$, $\CM_{(r,s)}(S)$ is the groupoid of triples
 $(A,\iota,\lambda)$
where  
 \begin{altitemize}
\item $(A,\iota)$ is an abelian scheme over $S$, with $O_{F'}$-action $\iota: O_{F'}\to \End(A)$ satisfying the Kottwitz condition of signature $(r, s)$, and 
\item $\lambda:A\to A^\vee$ is a polarization whose Rosati involution induces on $O_{F'}$ the non-trivial Galois automorphism of $F'/F$, and such that $\ker(\lambda)$ is contained in $ A[\iota(\varpi)]$ of rank $\#(O_{F'}/(\varpi))^{n}$ (resp.\ $\#(O_{F'}/(\varpi))^{n-1}$) when $n=r+s$ is even (resp.\ odd).
In particular, we have $\ker(\lambda)=A[\iota(\varpi)]$ if $n=r+s$ is even.
\end{altitemize}
 Now we assume that $(r,s)=(1,n-1)$ or $(n-1,1)$. We further impose the {\em wedge condition} and the {\em (refined) spin condition}, cf.\ \cite[\S4.4]{RSZ3}. The functor is represented by a Deligne--Mumford stack again denoted by $\CM_{(r,s)}$.
It is {\em smooth} over $\Spec O_{F'}$, despite the ramification of the field extension $F'/\BQ$, cf.\ \cite[\S4.4]{RSZ3}. 
Then we have an integral model of copies of the Shimura variety $\Sh_{ K_{\wt \RG}^\circ}(\wt \RG)$ defined by
$$
\CM_{ K_{\wt \RG}^\circ}(\wt \RG)=\CM_{(0,1)} \times_{\Spec O_{F'}} \CM_{(1,n-2)}\times_{\Spec O_{F'}}\CM_{(1,n-1)}.
$$  (In \cite[\S5.1]{RSZ3} we do cut out the desired Shimura variety with the help of a sign invariant. Here, implicitly we need to replace this space by its toroidal compactification. )

We now describe the arithmetic diagonal cycle (or rather, its integral model) for the level $K_{\wt \RH}^\circ=K_{\wt \RG}^\circ\cap \wt \RH(\BA_f)$. When $n$ is odd (so $n-1$ is even), we define 
$$
 \CM_{K_{\wt \RH}^\circ} \bigl(\wt \RH\bigr)=\CM_{(0,1)} \times_{\Spec O_{F'}} \CM_{(1,n-2)},
 $$ and we can define an embedding explicitly by ``taking products" (one sees easily that the conditions on the kernels of polarizations are satisfied):
\begin{equation}\label{inc M}
   \xymatrix@R=0ex{
	   \CM_{K_{\wt \RH}^\circ} \bigl(\wt \RH\bigr) \ar[r]  & \CM_{K_{\wt\RG}^\circ}\bigl(\wt \RG\bigr)\\
		\bigl( A_0,\iota_0,\lambda_0,A^\flat,\iota^\flat,\lambda^\flat \bigr) \ar@{|->}[r]
		   & \bigl( A_0,\iota_0,\lambda_0,A^\flat,\iota^\flat,\lambda^\flat, A^\flat \times  A_0,\iota^\flat \times \iota_0,\lambda^\flat \times \lambda_0\bigr) .
	}
\end{equation}	
When  $n$ is even, the situation is more subtle; see \cite[\S4.4]{RSZ3}.

 With the smooth integral model, we have an unconditionally defined height pairing \eqref{eqn BB} on $X=\Sh_{ K_{\wt \RG}^\circ}(\wt \RG)$. Now we again apply a suitable Hecke--Kunneth projector to the cycle $z_K$ for $K=K_{\wt \RG}^\circ$, and we obtain a cohomologically trivial cycle $z_{K,0}\in\Ch(\Sh_{ K_{\wt \RG}^\circ}(\wt \RG))_0 $. We define
  \begin{align}\label{int f}
 \Int(f)=\Big( R(f)\ast z_{K,0},\,\,z_{K,0}\Big)_{{\rm BB}},\quad f\in \sH\left(\wt\RG, K_{\wt \RG}^\circ\right),
 \end{align}
 where $R(f)$ is the associated Hecke correspondence. Let $\sH^{{\rm ram}_{F'}}(\wt \RG) $ be the spherical Hecke algebra away from the set ${\rm ram}_{F'}$ of primes ramified in $F'/F$.

 Parallel to Conjecture \ref{conj dis}, we can state an alternative version of the arithmetic Gan--Gross--Prasad conjecture in terms of spherical characters  for the special level $K_{\wt \RG}^\circ$.
\begin{conj}\label{conj:sph AGGP}
There is a decomposition $$
 \Int(f)=\sum_{\pi}\Int_\pi(f), \quad\text{for all }\, f\in \sH\left(\wt\RG, K_{\wt \RG}^\circ\right),
 $$where the sum runs over all automorphic representations of $\wt\RG(\BA)$ that appear in the cohomology $H^{\ast}(\Sh_{ \wt G, K_\RG})$ and are  trivial on $Z^\BQ(\BA)$, and $\Int_\pi$ is an eigen-distribution for the spherical Hecke algebra $\sH^{{\rm ram}_{F'}}(\wt \RG) $ with eigen-character $\lambda_{\pi^{ {\rm ram}_{F'} }}$ in the sense of \eqref{eigen dis}.
 
If such a representation $\pi$ is tempered, then 
 $$
 \Int_\pi(f)=2^{-\beta_\pi} \sL'(1/2,\pi) \prod_{v<\infty}\BI_{\pi_v}(f_v).
 $$
\end{conj}
Here the constant $\beta_\pi$ is the same as  in  \S\ref{sss:sph GGP}, and there is a natural extension of the local spherical characters $\BI_{\pi_v}$ to the triple $(\wt \RG,\wt \RH,\wt \RH)$.

 \begin{remark}Conjecture \ref{conj:sph AGGP} can be viewed as a refined version of the arithmetic Gan--Gross--Prasad conjecture \ref{conj AGGP}. Other refinements were also given in by the author \cite{Z09} and independently by S. Zhang \cite{ZSW10}. Both of them rely on standard conjectures on height pairings, and hence are conditional. 
 \end{remark}

\section{Shtukas and higher Gross--Zagier formula}
\label{s:HGZ}
Now let $F$ be the function field of a geometrically connected smooth proper curve $X$ over a finite field $k=\BF_q$. We may consider the analog of the special pair of Shimura data (cf.\ \ref{sss:pair Sh}) in the context of Shtukas. Now there is much more freedom since we do not have the restriction from the archimedean place. One may choose an $r$-tuple of coweights of $\RG$ to define $\RG$--Shtukas (with $r$-modifications), and the resulting moduli space lives over the $r$-fold power 
$$
X^r=\underbrace{X\times_{\Spec k}\dotsc\times_{\Spec k} X }_{r \text{ times}}.
$$
 This feature is completely missing in the number field case, where we only have two available options:
	\begin{altenumerate}
	\item when $r=0$, the automorphic quotient $[\RG]_K$ plays an analogous role, cf.\ \ref{sss:aut quot}.
	\item when $r=1$, we have Shimura varieties $\Sh_{K_{\RG}}(\RG)$ associated to a Shimura datum $\left(\RG,\{h_\RG\}\right)$. These varieties live over $\Spec E$ for a number field $E$.
\end{altenumerate}

 In \cite{YZ1} Yun and the author studied a simplest case, i.e., the special cycle on the moduli stack of rank two Shtukas with arbitrary number $r$ of modifications.  We connect their intersection numbers to the $r$-th order derivative of certain L-functions. We may view the result as an analog of Waldspurger's formula (for $r=0$) and the Gross--Zagier formula (for $r=1$).

\subsection{The Heegner--Drinfeld cycle}

Let $\RG=\PGL_2$ and let $\Bun_2$ be the stack of rank two vector bundles on $X$. The Picard stack $\Pic_X$ acts on $\Bun_2$ by tensoring the line bundle. Then $\Bun_\RG=\Bun_2/\Pic_X$ is the moduli stack of $\RG$-torsors over $X$.

Let $r$ be an {\em even} integer.  Let $\mu\in\{\pm\}^{r}$ be an $r$-tuple of signs such that exactly half of them are equal to $+$. Let  $\Hk_2^\mu$  be the Hecke stack, i.e., $\Hk_2^\mu(S)$ is the groupoid of data 
$$
(\CE_0,\cdots,\CE_r,x_1,\cdots,x_r, f_1,\cdots,f_r),
$$
 where the $\CE_i$'s are vector bundles of rank two over $X\times S$, the $x_i$'s are $S$-points of $X$, each $f_i$ is a minimal upper (i.e., increasing) modification if $\mu_i=+$, and minimal lower (i.e., decreasing) modification if $\mu_i=-$, and the $i$-th modification takes place along the graph of $x_i:S\to X$,
\begin{equation*}
\xymatrix{\calE_{0}\ar@{-->}[r]^{f_{1}} &\calE_{1}\ar@{-->}[r]^{f_{2}}& \cdots \ar@{-->}[r]^{f_{r}}&\calE_{r}}.
\end{equation*}
The Picard stack $\Pic_X$ acts on $\Hk_2^\mu$ by simultaneously tensoring the line bundle. 
Define $\Hk^\mu_\RG=\Hk_2^\mu/\Pic_X$. Recording $\CE_i$ defines a projection  $p_i: \Hk^\mu_\RG\to \Bun_\RG$. 

The moduli stack $\Sht^{\mu}_\RG$ of Drinfeld $\RG$-Shtukas of type $\mu$ for the group $\RG$ is defined by the cartesian diagram
\begin{equation}\label{ShtrG}
\xymatrix{\Sht^{\mu}_{\RG}\ar[d]\ar[rr] && \Hk^{\mu}_{\RG}\ar[d]^{(p_{0},p_{r})}\\
\Bun_{\RG}\ar[rr]^{(\id,\Fr)} && \Bun_{\RG}\times \Bun_{\RG}}.
\end{equation}
The stack $\Sht^{\mu}_{\RG}$ is a  Deligne--Mumford stack over $X^r$, and the natural morphism 
$$
\xymatrix{\pi_{\RG}^\mu: \Sht_{\RG}^\mu\ar[r] &X^r}$$ 
is smooth of relative dimension $r$.

Let
$\nu: X'\to X$ be a finite \'etale cover of degree $2$ such that $X'$ is also geometrically connected.  Denote by $F'=k(X')$ the function field.  Let $\RT=(\RR_{F'/F}\BG_m)/\BG_m$ be the non-split torus associated to the double cover $X'$ of $X$. The stack $\Sht^{\mu}_{\RT}$ of $\RT$-Shtukas is defined analogously, with the rank two bundles $\CE_{i}$ replaced by line bundles $\CL_i$ on $X'$, and the points $x_{i}$ on $X'$. 
Then we have a map
$$
\xymatrix{\pi^{\mu}_{\RT}: \Sht^{\mu}_{\RT}\ar[r]& X'^r}$$
which is a torsor under the finite Picard stack $\Pic_{X'}(k)/\Pic_{X}(k)$. In particular,  $\Sht^{\mu}_{\RT}$ is a proper smooth Deligne--Mumford stack over $\Spec k$.

There is a natural finite morphism of stacks over $X^r$, induced by the natural map $\nu_\ast:\Pic_{X'}\to \Bun_2$,
$$
\xymatrix{\Sht^{\mu}_{\RT}\ar[r]& \Sht^{\mu}_{\RG}}.
$$
It induces a finite morphism
$$
\xymatrix{\theta^{\mu}:\Sht^{\mu}_{\RT}\ar[r]& \Sht'^{\mu}_{\RG}:=\Sht_{\RG}^{\mu}\times_{X^r}X'^r}.
$$
This defines a class in the Chow group of proper cycles of dimension $r$ with $\BQ$-coefficients,
\begin{align}\label{HD cyc}
Z_{\RT}^{\mu}:=\theta^{\mu}_{*}[\Sht^{\mu}_{\RT}]\in \Ch_{c,r}(\Sht'^{\mu}_{\RG})_\BQ.
\end{align}
In analogy to the Heegner cycle in the Gross--Zagier formula \cite{GZ,YZZ1} in the number field case, we call $Z_{\RT}^{\mu}$ the {\em Heegner--Drinfeld cycle} in our setting.

\begin{remark}
The construction of the Heegner--Drinfeld cycle extends naturally to higher rank Shtukas (of rank $n$ over $X'$, respectively rank $2n$ over $X$) of type $\mu=(\mu_1,\cdots, \mu_r)$. Here $\mu_i$ are coweights of $\GL_{n}$ (or $\GL_{2n}$) given by $(\pm1,0,\cdots,0)\in \BZ^n$ (or $\BZ^{2n}$). 
\end{remark}

\subsection{Taylor expansion of $L$-functions} 
 Consider the middle degree cohomology with compact support 
$$V=\cohoc{2r}{(\Sht'^{\mu}_{G})\otimes_{k}\kbar,\ov\BQ_\ell}(r).$$
 This vector space is endowed with the cup product 
$$
\xymatrix{(\cdot,\cdot): V\times V\ar[r]&\ov\BQ_\ell}.
$$

Let $\ell$ be a prime number different from $p$.   Let $K=\prod_{v} \RG(O_{F_v})$, 
  and let $\sH_{\ov\BQ_\ell}=\sH(\RG(\BA),K)$ be the spherical Hecke algebra with $\ov\BQ_\ell$-coefficients. For any maximal ideal $\fkm\subset\sH_{\ov\BQ_\ell}$, we define the generalized eigenspace of $V$ with respect to $\fkm$ by
\begin{equation*}
V_{\fkm}=\cup_{i>0}V[\fkm^i].
\end{equation*}
We also define a subspace $V_{\Eis}$ with the help of an {\em Eisenstein ideal}, cf \cite[\S4.1.2]{YZ1}. 
Then we prove that there is a spectral decomposition, i.e., an orthogonal decomposition of $\sH_{\ov\BQ_\ell}$-modules,
\begin{equation}\label{intro V decomp}
V=V_{\Eis}\oplus\left(\bigoplus_{\fkm}V_{\fkm}\right),
\end{equation}
where $\fkm$ runs over a finite set of maximal ideals of $\sH_{\ov\BQ_\ell}$, and each $V_{\fkm}$ is an $\sH_{\ov\BQ_\ell}$-module of finite dimension over $\ov\BQ_\ell$ supported at the maximal ideal $\fkm$; see \ \cite[Thm. 7.16]{YZ1} for a more precise statement.

Let $\pi$ be an everywhere unramified cuspidal automorphic representation $\pi$ of $\RG(\BA)$. The standard $L$-function $L(\pi,s)$ is a polynomial of degree $4(g-1)$ in $q^{-s-1/2}$, where $g$ is the genus of $X$. Let $\pi_{F'}$ be the base change to $F'$, and let $L(\pi_{F'},s)$ be its standard $L$-function. Let $L(\pi,\Ad,s)$ be the adjoint $L$-function of $\pi$ and define
\begin{equation}
\label{eqn sL}
\sL(\pi_{F'},s)= \epsilon(\pi_{F'},s)^{-1/2}\frac{L(\pi_{F'},s)}{L(\pi,\Ad,1)},
\end{equation}
where the the square root is understood as
$ \epsilon(\pi_{F'},s)^{-1/2}=q^{4(g-1)(s-1/2)}.$
In particular, we have a functional equation:
$$
\sL(\pi_{F'},s)=\sL(\pi_{F'},1-s).
$$
We consider the Taylor expansion at the central point $s=1/2$:
$$
\sL(\pi_{F'},s)= \sum_{r\geq 0} \sL^{(r)}(\pi_{F'},1/2)\frac{(s-1/2)^r}{r!},
$$
i.e.,
$$
\sL^{(r)}(\pi_{F'},1/2)=\frac{d^r}{ds^r}\Big|_{s=0}  \left(   \epsilon(\pi_{F'},s)^{-1/2}\frac{L(\pi_{F'},s)}{L(\pi,\Ad,1)}\right).
$$
If $r$ is odd, by the functional equation we have
$$
\sL^{(r)}(\pi_{F'},1/2)=0.
$$

Now we fix an isomorphism $\BC\simeq \ov\BQ_\ell$. Let $\fkm=\fkm_\pi$ be the kernel of the associated character $\l_{\pi}:\sH_{\ov\BQ_\ell}\to\ov\BQ_\ell$, and rename $V_\fkm$ in  \eqref{intro V decomp} as $V_\pi$. Then our main result in \cite{YZ1} relates the $r$-th Taylor coefficient to the self-intersection number of the $\pi$-component of the Heegner--Drinfeld cycle  $\theta^{\mu}_{*}[\Sht^{\mu}_{\RT}]$. 
\begin{thm}\label{th:main cycle} Let $\pi$ be an everywhere unramified cuspidal automorphic representation of $\RG(\BA)$.
Let $[\Sht^{\mu}_{\RT}]_{\pi}\in V_{\pi}$ be the projection of the cycle class of  $\cl\left(\theta^{\mu}_{*}[\Sht^{\mu}_{\RT}]\right)\in V$ to the direct summand $V_{\pi}$ under the decomposition \eqref{intro V decomp}. Then
$$
 \frac{1}{2(\log q)^{r}} \,\lvert \omega_X\rvert \,\sL^{(r)}\left(\pi_{F'},1/2\right) =\Big([\Sht^{\mu}_{\RT}]_{\pi},\quad [\Sht^{\mu}_{T}]_{\pi} \Big),
$$
where $\omega_X$ is the canonical divisor, and $\lvert \omega_X\rvert=q^{-2g+2}$.
\end{thm}

\begin{remark} Here we only consider \'etale double covers $X'/X$, and everywhere unramified $\pi$ (whence the $L$-function has nonzero Taylor coefficients in even degrees only). In \cite{YZ2}, Yun and the author are extending the theorem above to the case when $X'/X$ is ramified at a finite set $R$ and $\pi$ has Iwahori levels at $\Sigma$ such that $R\cap \Sigma=\emptyset$.
\end{remark}  

\subsection{Comparison with the conjecture of Birch and Swinnerton-Dyer}
Let $\pi$ be as in Theorem \ref{th:main cycle}, and $\rho_\pi$ the associated local system of rank two over the curve $X$ by the global Langlands correspondence. Let $
W'_\pi= H^1\left( X'\times \ov k,\, \rho_\pi\right),
$ a $\ov\BQ_\ell$-vector space with the Frobenius endomorphism $\Fr$. The L-function $L(\pi_{F'},s)$ is then given by
$$
L(\pi_{F'},s-1/2)=\det\left(1-q^{-s}\Fr\bigm| W_\pi'\right).
$$In particular, the dimension of the eigenspace $W'^{\Fr=q}_\pi$ is at most $\ord_{s=1/2}L(\pi_{F'},s)$ (the conjectural semi-simplicity of $\Fr$ implies an equality).
 It is {\em expected} that, the complex ${\bf R}\pi^{\mu}_{\RG,!}\ov\BQ_\ell$ on $X^r$ decomposes as a direct sum of $\sH_{\ov\BQ_\ell}$-modules
$$
{\bf R}\pi^{\mu}_{\RG,!}\ov\BQ_\ell=\Big(\bigoplus_{\pi\text{\,cuspidal}} \pi^K\otimes\big(\underbrace{\rho_\pi \boxtimes\cdots\boxtimes \rho_\pi}_{r\text{\,\quad times}}\big)\Big) \bigoplus\text{ ``a direct summand"},
$$
such that $V_\pi=V_{\fkm_\pi}$ in \eqref{intro V decomp} (for $\fkm_\pi=\ker(\lambda_\pi)$) corresponds to $\pi^K\otimes W'^{\otimes r}_\pi.$ 
From now on we assume this decomposition.
Then the cohomology class of the Heegner--Drinfeld cycle defines an element $Z_\pi^\mu\in \pi^K\otimes W'^{\otimes r}_\pi$. One can show that $Z_\pi^\mu$ is an eigen-vector for the operator $\id\otimes \Fr^{\otimes r}$ with eigenvalue $q^r$.
Then Theorem \ref{th:main cycle} shows that this class does not vanish when $r\geq \ord _{s=1/2}L(\pi_{F'},s)$, provided that $L(\pi_{F'},s)$ is not a constant. \begin{conj}\label{conj basis}
Let $r=\ord _{s=1/2}L(\pi_{F'},s)$. The class $Z_\pi^\mu$ belongs to $\pi^K\otimes \wedge^r\left(W'^{\Fr=q}_\pi\right)$.
\end{conj} Note that the generalization of the conjecture of Birch and Swinnerton-Dyer to function fields by Artin and Tate predicts that 
$
\dim W'^{\Fr=q}_\pi=\ord _{s=1/2}L(\pi_{F'},s).
$

We have a similar conjecture when $X'/X$ is ramified at a finite set $R$ and $\pi$ has Iwahori levels at $\Sigma$ such that $R\cap \Sigma=\emptyset$, cf.\ \cite{YZ2}.
In a forthcoming work, Yun and the author plan to prove that
\begin{altenumerate}
\item Let $r_0\geq 0$ be the smallest integer $r$ such that $Z_\pi^\mu\neq 0$ for some $\mu\in\{\pm\}^{r}$. Then $\dim W'^{\Fr=q}_\pi=r_0$, and the class $Z_\pi^\mu$ gives a  basis to the line $\pi^K\otimes \wedge^r\left(W'^{\Fr=q}_\pi\right)$.
\item {\em $\ord _{s=1/2}L(\pi_{F'},s)=1$ if and only if $\dim W'^{\Fr=q}_\pi=1$}. In particular, if $\ord _{s=1/2}L(\pi_{F'},s)=3$, then $\dim W'^{\Fr=q}_\pi=3$.
\end{altenumerate}

\section{Relative trace formula}
\subsection{An overview of RTF}

A natural tool to study automorphic period integrals  is the {\em relative trace formula} (RTF) introduced by Jacquet. For the reader's convenience, we give a very brief overview of the relative trace formula (cf.\ the survey articles \cite{J05,La06,La10,Of09}). 
 
 We start with a triple $(\RG,\RH_1,\RH_2)$ consisting of a reductive group $\RG$ and two subgroups $\RH_1,\RH_2$ defined over $F$. Known examples suggest that we may further assume that the pairs $(\RH_i,\RG)$ are spherical (cf.\ \S\ref{sss:sph}), although this is not essential to our informal discussion here.
  
 To a test function $f\in \mathscr{C}_c^\infty(\RG(\BA))$ we associate an automorphic kernel function,
$$
K_{f}(x,y):=\sum_{\gamma\in \RG(F)}f(x^{-1}\gamma y),\quad x,y\in \RG(\BA),
$$
which is invariant under $\RG(F)$ for both variables $x$ and $y$. This defines an integral operator representing $R(f)$ for the action $R$ of $\RG(\BA)$ on the Hilbert space $L^2([\RG])$. Therefore the kernel function has a spectral decomposition, and the contribution of a cuspidal automorphic representation $\pi$ to the kernel function is given by
\begin{equation}
\label{eqn K pi}
K_{\pi,f}(x,y)=\sum_{\varphi\in \mathrm{OB}(\pi)}\bigl(\pi(f)\varphi \bigr)(x)\,\ov{ \varphi(y)},
\end{equation}
where the sum runs over an orthonormal basis $\mathrm{OB}(\pi)$ of $\pi$ (with respect to the Petersson inner product). 

Then we consider a linear functional on $\mathscr{C}_c^\infty(\RG(\BA))$,
\begin{equation}
\label{eqn I(f)}
\BI(f)=\int_{[\RH_1]}\int_{[\RH_2]}K_f(h_1,h_2)\,dh_1\,dh_2.
\end{equation}
The spectral contribution \eqref{eqn K pi} from an automorphic representation $\pi$ is the (global) {\it spherical character} (relative to $(\RH_1,\RH_2)$), denoted by $\BI_\pi(f)$. Similar to \eqref{global char}, this is equal to 
$$
\BI_\pi(f)=\sum_{\phi\in \mathrm{OB}(\pi)}\sP_{\RH_1}(\pi(f)\phi)\ov{\sP_{\RH_2}(\phi)}.
$$

 Let $\RH_{1,2}:= \RH_1\times \RH_2.$
Then $\RH_{1, 2}$ acts on $\RG$ by $(h_1, h_2)\colon \gamma\mapsto h_1^{-1}\gamma h_2$. For certain nice orbits $\gamma\in \RG(F)/\RH_{1,2}(F)$, we can define orbital integrals (relative to $\RH_{1,2}$):
\begin{align}\label{orb}
\Orb(\gamma,f):=\vol([\RH_{1,2,\gamma}])\int_{ \RH_{1,2}(\BA)/\RH_{1,2,\gamma}(\BA)}f(h_1^{-1}\gamma h_2)\,dh_1\,dh_2,
\end{align}
where $\RH_{1,2,\gamma}$ denotes the stabilizer of $\gamma$, and $\vol$ stands for ``volume". 

 The relative trace formula attached to the triple $(\RG,\RH_1,\RH_2)$ is then the identity between the spectral expansion and the geometric expansion of $\BI(f)$:
 $$
 \sum_{\pi}\BI_\pi(f)+\cdots=\sum_{\gamma}\Orb(\gamma,f)+\cdots,
 $$
 where the $\cdots$ parts need more care (in fact saying so is oversimplifying). 
We will use  ${\rm RTF}_{(\RG,\RH_1,\RH_2)}$ to stand for the above relative trace formula identity. This is only a very coarse form, and depending on the triple $(\RG,\RH_1,\RH_2)$ the identity may need further refinements such as stabilization, as experience with the Arthur--Selberg trace formula suggests.

\begin{remark}When  we take the triple $(\RH\times \RH,\Delta_\RH,\Delta_\RH),$
where $\Delta_\RH\subset \RH\times\RH$ is the diagonal embedding of $\RH$, the associated relative trace formula is essentially equivalent to the Arthur--Selberg trace formula associated to $\RH$. Therefore the RTF can be viewed as a generalization of the Arthur--Selberg trace formula. 
\end{remark}

In application to questions such as the Gan--Gross--Prasad conjecture, we need to compare two RTFs that are close to each other,
$$
{\rm RTF}_{(\RG,\RH_1,\RH_2)}\longleftrightarrow {\rm RTF}_{(\RG',\RH'_1,\RH'_2)}.
$$
The comparison allows us to connect the automorphic periods on $\RG$ to those on $\RG'$.
There are many successful examples, although it is a subtle question how to seek comparable RTFs in general. 

\subsection{Jacquet--Rallis RTFs}
\label{ss:JR}
\subsubsection{}
We recall the two RTFs constructed by Jacquet and Rallis \cite{JR} to attack the Gan--Gross--Prasad conjecture in the Hermitian case (cf.\ \S\ref{sss:GGP}).  

The first RTF deals with the automorphic $\RH$-period integral on $\RG$ and is associated to the triple $(\RG,\RH,\RH)$. The second one is associated to the triple $(\RG',\RH'_1,\RH'_2)$ where
$$
\RG'=\RR_{F'/F}(\GL_{n-1}\times \GL_{n}),$$
and $$\RH'_1=\RR_{F'/F}\GL_{n-1},\quad \RH'_2=\GL_{n-1}\times\GL_{n},
$$
where $(\RH'_1,\RG')$ is the Rankin--Selberg pair, and $(\RH'_2,\RG')$ the Flicker--Rallis pair, cf.\ \S\ref{sss:sph}.
Moreover it is necessary to insert a quadratic character of $\RH_2'(\BA)$:
$$
\eta=\eta_{n-1,n}:(h_{n-1},h_{n
})\in \RH'_2(\BA)\mapsto \eta^{n-2}_{F'/F}(\det(h_{n-1}))\eta^{n-1}_{F'/F}(\det(h_{n})),
$$
where $\eta_{F'/F}: F^\times\bs \BA^\times\to\{\pm 1\}$ is the quadratic character associated to $F'/F$ by class field theory. 

For the later application to the arithmetic Gan--Gross--Prasad conjecture, we introduce (cf.\  \cite[\S3.1]{Z12}) the global distribution on $\RG'(\BA)$ parameterized by a complex variable $s\in\BC$,
\begin{align}
\label{J f'}
\BJ(f',s)=\int_{[\RH'_1]}\int_{[\RH_2']}K_{f'}(h_1,h_2)\,\bigl|\det(h_1)\bigr|^s\,\eta(h_2)\,dh_1\,dh_2,\quad f'\in \mathscr{C}^\infty_c(\RG'(\BA)).
\end{align}
We set
\begin{align*}
\BJ(f')=\BJ(f',0).
\end{align*}

\begin{remark}Due to the presence of the Flicker--Rallis pair $(\RH'_2,\RG')$, the cuspidal part of the spectral side in ${\rm RTF}_{(\RG',\RH'_1,\RH'_2)}$ only contains those  automorphic representations that are in the image of the quadratic base change from unitary groups. This gives the hope that the spectral sides of the two RTFs should match.
\end{remark}

\begin{remark}
In \cite{GGP}, Gan, Gross, and Prasad also made global conjectures for  $\SO_n\times\SO_m$ and  $\U_n\times\U_m$ when $|n-m|>1$. Towards them in the Hermitian cases, Y. Liu in \cite{Liu1,Liu2} has generalized the construction of Jacquet and Rallis.  
\end{remark}

\subsubsection{Geometric terms: orbital integrals}\label{sss:JR orb}

In the comparison of geometric sides of two RTFs, we need to match orbits and orbital integrals.  We review the comparison in the Jacquet--Rallis case. 

We call an element $\gamma\in \RG(F)$ \emph{regular semi-simple} (relative to the action of $\RH_{1, 2}=\RH_1\times\RH_2$) if its orbit under $\RH_{1, 2}$ is Zariski closed and its stabilizer is of minimal dimension. 
The regular semi-simple orbits will be the nice ones for the study of orbital integrals. For the triples $(\RG,\RH,\RH)$ and $(\RG',\RH'_1,\RH'_2)$ in the Jacquet--Rallis RTFs, the condition is equivalent to $\gamma$ having Zariski closed orbit and trivial stabilizer. In particular, for such $\gamma$ the orbital integral \eqref{orb} simplifies. We denote by $\RG(F)_\rs$ (resp.\ $\RG'(F)_\rs$) the set of regular semi-simple elements in $\RG(F)$ (resp.\ $\RG'(F)$). We denote by $\bigl[\RG(F)_\rs\bigr]$ and $\bigl[\RG'(F)_\rs\bigr]$ the respective sets of orbits. 

 Depending on the pair of Hermitian spaces $W:=(W_{n-1},W_n)$, we denote the triple $(\RG,\RH,\RH)$ by $(\RG_W,\RH_W,\RH_W)$. We consider the equivalence relation $(W_{n-1},W_n)\sim (W_{n-1}',W_n')$ if there is a scalar $\lambda\in F^\times$ such that we have isometries $W_{n-1}'\simeq {}^\lambda W_{n-1}$ and $W_{n}'\simeq \,^\lambda W_{n}$, where the left superscript changes the Hermitian form by a multiple $\lambda$. 
There is a natural bijection (cf.\ \cite[\S2]{Z12}, \cite[\S2]{RSZ1}) 
\begin{align}\label{orb bij}
   \coprod_{W}\,\, \bigl[\RG_{W}(F)_\rs\bigr]\isoarrow  \bigl[\RG'(F)_\rs\bigr],
\end{align}
where the left hand side runs over all pairs $W=(W_{n-1},W_n)$ up to equivalence.
This bijection holds for any quadratic extension of fields $F'/F$ of characteristic not equal to $2$.

Now we let $F'/F$ be a quadratic extension of {\em local fields}. For $g\in \RG_W(F)_\rs$ and $f\in \mathscr{C}^\infty_c(\RG_W(F))$ we introduce the orbital integral
\begin{align}\label{orb U}
\Orb(g,f)=
\int_{ \RH_{1,2}(F)}f(h_1^{-1}g h_2)\,dh_1\,dh_2.
\end{align}
For $\gamma\in \RG'(F)_\rs$, $f'\in \mathscr{C}^\infty_c(\RG'(F))$, and $s\in\BC$, we introduce the (weighted) orbital integral
\begin{align}\label{orb GL}
\Orb(\gamma,f',s)=
\int_{ \RH'_{1,2}(F)}f(h_1^{-1}\gamma h_2)\bigl|\det(h_1)\bigr|^s\eta(h_2)\,dh_1\,dh_2.
\end{align}
We set
\begin{equation}\label{orb del}
   \Orb(\gamma,f') := \Orb(\gamma,f', 0)  \quad\text{and}\quad \del(\gamma,f') := \frac d{ds} \Big|_{s=0} \Orb(\gamma, f',s)  .
\end{equation}

\begin{defn}\label{def tran}
\begin{altenumerate}
\item
We say that $g\in \RG_W(F)_\rs$ and $\gamma\in \RG'(F)_\rs$ {\it match} if their orbits correspond to each other under \eqref{orb bij}. 
\item Dually, we say that a function $f'\in \mathscr{C}^\infty_c(\RG'(F))$ and a tuple $\{f_W\in \mathscr{C}^\infty_c(\RG_W(F))\}$ indexed by $W$ (up to equivalence) are {\em transfers} of each other if  for each $W$ and each $g\in \RG_W(F)_\rs$,
$$
 \Orb(g,f_W)=\omega(\gamma)  \Orb(\gamma,f')
$$
whenever $\gamma\in \RG'(F)_\rs$ matches $g$. Here $\omega(\gamma) $ is a certain explicit transfer factor \cite{Z14a,RSZ2}. 
\item We say that a component $f_W$ in the tuple is a {\it transfer} of $f'$ if the remaining components of the tuple are all zero.
\end{altenumerate}
\end{defn}

In \cite{Z14a} we prove the following.
\begin{thm}\label{thm ST}
Let $F'/F$ be a quadratic extension of $p$-adic local fields (then there are two equivalence classes of pairs of Hermitian spaces denoted by $W,W^\flat$ respectively). Then for any $f'\in \mathscr{C}^\infty_c(\RG'(F))$ there exists a transfer  $(f_0,f_1)\in \mathscr{C}^\infty_c(\RG_{W}(F))\times \mathscr{C}^\infty_c(\RG_{W^\flat}(F))$, and for any pair $(f_0,f_1)\in \mathscr{C}^\infty_c(\RG_{W}(F))\times \mathscr{C}^\infty_c(\RG_{W^\flat}(F))$ there exists a transfer $f'\in \mathscr{C}^\infty_c(\RG'(F))$. 
\end{thm}
This was conjectured by Jacquet and Rallis \cite{JR}. For archimedean local fields $F'/F$, an ``approximate transfer" is proved by Xue \cite{Xue}.

\subsubsection{Spectral terms: spherical characters}\label{sss:JR sph}
We are now back to $F$ being a number field.
For the triple $(\RG,\RH_1,\RH_2)$ we have defined the global (resp.\ local) spherical characters by \eqref{global char} (resp.\ by \eqref{def local char G}). For the triple $(\RG',\RH'_1,\RH'_2)$, we define the global spherical character associated to a cuspidal automorphic representation $\Pi$ of $\RG'(\BA)$,
$$
\BJ_\Pi(f',s)=\sum_{\phi\in \mathrm{OB}(\Pi)}\sP_{\RH'_1,s}(\Pi(f)\phi)\ov{\sP_{\RH'_2,\eta}(\phi)}, \quad f'\in \sC_c^\infty(\RG'(\BA)),\quad s\in\BC,
$$
where $\sP_{\RH'_1,s}$ is the automorphic period integral  $\sP_{\RH'_1,\chi_s}$ for the character $\chi_s$ of $\RH_1(\BA)\in\GL_{n-1}(\BA_{F'})$ defined by  $h\mapsto  |\det(h)|_{F'}^s$. We set $\BJ_{\Pi}(f')=\BJ_{\Pi}(f',0)$. We expect to have a global character identity \cite[Conj. 4.2]{Z14b}:
\begin{conj}\label{conj global char}
Let $\pi$ be a tempered cuspidal automorphic representation of $\RG(\BA)$ such that $\Hom_{\RH(\BA)}(\pi,\BC)\neq 0$. Let $\Pi={\rm BC}(\pi)$ be its base change which we assume is cuspidal. Then for $f'\in \mathscr{C}^\infty_c(\RG'(\BA))$ and any transfer $f\in \mathscr{C}^\infty_c(\RG(\BA))$, 
$$
\BI_{\pi}(f)=2^{-\beta_\pi}\BJ_{\Pi}(f').
$$
\end{conj}
\begin{remark}\label{rem spec}
In fact here we have $\beta_\pi=2$ due to the cuspidality of $\Pi={\rm BC}(\pi)$.  Conjecture \ref{conj global char} is known when $\pi_v$ is supercuspidal at a place $v$ split in $F'/F$, cf.\ \cite{Z14a}. In general it should follow from a full spectral decomposition of the Jacquet--Rallis relative trace formulas, and perhaps along the way one will discover the correct definition of the global spherical character $\BJ_{\Pi}(f')$ when $\Pi={\rm BC}(\pi)$ is not cuspidal.
\end{remark}

In \cite[\S3.4]{Z14b}, we defined a local spherical character $\BJ_{\Pi_v}(f'_v,s)$  for any tempered $\Pi_v$ (depending on some auxiliary data). Then,  for pure tensors $f'=\otimes_v f'_v$, the global spherical character decomposes naturally as an Euler product for $\Pi={\rm BC}(\pi)$,
\begin{align}\label{char glo2loc}
\BJ_\Pi(f',s)=2^{-\beta_\pi}\sL(s+1/2,\pi)\prod_{v}\BJ_{\Pi_v}(f'_v,s).
\end{align}
We set $\BJ_{\Pi_v}(f'_v)=\BJ_{\Pi_v}(f'_v,0)$. We expect to have a local character identity \cite[Conj. 4.4]{Z14b}:
\begin{conj}\label{conj local char}
Let $\pi_v$ be a tempered representation of $\RG(F_v)$ such that $\Hom_{\RH(F_v)}(\pi,\BC)\neq 0$. Let $\Pi_v={\rm BC}(\pi_v)$ be its base change. Then for $f'_v\in \mathscr{C}^\infty_c(\RG'(F_v))$ and any transfer $f_v\in \mathscr{C}^\infty_c(\RG(F_v))$,
$$\BJ_{\Pi_v}(f'_v)
=\kappa_{\pi_v}\,\BI_{\pi_v}(f_v),
$$
where $\kappa_{\pi_v}$ is an explicit constant.
\end{conj}

By a theorem of Harish-Chandra, the character of an admissible representation of $\RG(F_v)$ for any reductive group $\RG$ over a $p$-adic local field $F_v$ admits a local expansion around the identity of $\RG(F_v)$ as a sum of Fourier transforms of unipotent orbital integrals. We have a partial analog for the local spherical character $\BJ_{\Pi_v}$.
\begin{thm}\label{thm:germ}Let $v$ be a non-archimedean place. Let $\Pi_v={\rm BC}(\pi_v)$ be the base change of a tempered representation $\pi_v$ of $\RG(F_v)$. For any small neighborhood of the identity element of $\RG'(F_v)$, there exists an admissible (in the sense of \cite[\S8.1]{Z14b}) function  $f'_v\in \mathscr{C}^\infty_c(\RG'(F_v))$ such that 
$$\BJ_{\Pi_v}(f'_v)=c_\Pi\cdot\wh{\mu}_{{\rm reg}}(f'_v),
$$
where $\wh{\mu}_{{\rm reg}}$ is the Fourier transform of the  (relative) regular unipotent orbital integral, cf.\ \cite[\S6.3, \S8.2]{Z14b}, and $c_\Pi$ is an explicit constant depending on $\Pi$. 
\end{thm}

We have the following.
\begin{thm}Conjecture \ref{conj local char} holds if $v$ is split in $F'/F$, or $F_v$ is a $p$-adic local field.
\end{thm}
The case of a split $v$ is rather easy \cite{Z14b}. The case of a supercuspidal representation $\pi_v$  was proved in \cite{Z14b}.  For the general $p$-adic case, the result is proved by Beuzart-Plessis in \cite{BP16}  using Theorem \ref{thm:germ}, a local relative trace formula for Lie algebras in \cite[\S4.1]{Z14a}, and a group analog in \cite{BP15b}. 

\begin{remark}
If $\Pi={\rm BC}(\pi)$ is cuspidal, then Conjecture \ref{conj global char} and Conjecture \ref{conj local char} together imply Conjecture \ref{conj dis}.
\end{remark}

\subsection{Arithmetic RTF}\label{ss:A RTF}
In  \cite{Z12}, the author introduced a relative trace formula approach to the arithmetic Gan--Gross--Prasad conjecture. Let $$
\partial\BJ(f')=\frac d{ds} \Big|_{s=0} \BJ(f',s),
$$cf.\ \eqref{J f'}.
Then the idea is that, in analogy to the usual comparison of two RTFs, we hope to compare the height pairing $\Int(f)$ in \eqref{int f} and  $\partial \BJ(f')$ for $f\in \sH(\wt\RG, K_{\wt \RG}^\circ)$  and ay transfer $f'\in \mathscr{C}^\infty_c(\RG'(\BA))$.
\begin{remark}
Here we note that there is no archimedean component in the test function $f$ on the unitary group side. Implicitly we demand that  $f'=\otimes_{v}f'_v\in \sC_c^\infty(\RG'(\BA))$, where $f'_\infty$ is a {\em Gaussian} test function in the sense of \cite[\S7.3]{RSZ3} (equivalently, we complete $f$ by tensoring a distinguished archimedean component $f_\infty$ as in \cite[\S3.2, (3.5)]{Z12}). We also note that, by the isomorphisms \eqref{proddec}, the orbits on $\wt \RG(F)_\rs$ are in natural bijection with those on $\RG(F)_\rs$, and all geometric terms related to $\RG$ in \S\ref{sss:JR orb} transport to $\wt \RG$. We will not repeat the definitions.
\end{remark}

To be able to work in a greater generality than the case the height pairing \eqref{eqn BB} is defined, in \cite{RSZ3}, Rapoport, Smithling, and the author turn to the arithmetic intersection theory $  \sform_\mathrm{GS}$ of 
Arakelov and Gillet--Soul\'e \cite[\S4.2.10]{GS} on the arithmetic Chow group $\wh\Ch^\ast(\CX)$ of a regular proper flat scheme (or DM stack) $\CX$ over $\Spec(O_{F'})$. For certain more general levels $K_{\wt \RG}\subset K_{\wt \RG}^\circ$, we construct regular integral models $\CM_{K_{\wt \RG}}(\wt \RG)$ of $\Sh_{K_{\wt \RG}}(\wt \RG)$ (essentially by adding Drinfeld level structures at split primes to the moduli space  $\CM_{K_{\wt \RG}^\circ}(\wt \RG)$, cf.\ \cite[\S4,\,\S5]{RSZ3}). We enhance the arithmetic diagonal cycle to an element $\wh z_{K}$ in $\wh\Ch^{n-1}(\CM_{K_{\wt \RG}}(\wt \RG))$, and we extend the action $R$ of a smaller Hecke algebra $\sH^{\rm spl}(\wt\RG, K_{\wt \RG})\subset \sH(\wt\RG, K_{\wt \RG})$ on the Chow group of $\Sh_{K_{\wt \RG}}(\wt \RG)$  to an action $\wh R$ on the arithmetic Chow group. We define 
\begin{align}
\Int(f)=\left(\wh R(f)\wh z_{K} , \,\wh z_{K}\right)_{{\rm GS}}, \quad f\in \sH^{\rm spl}(\wt\RG, K_{\wt \RG}).
\end{align}
We can then state an arithmetic intersection conjecture for the arithmetic diagonal cycle on the global integral model $\CM_{K_{\wt \RG}}(\wt \RG)$ \cite[\S8.1,\,\S8.2]{RSZ3}.
\begin{conj}\label{conj I=dJ}
Let $f\in \sH^{\rm spl}(\wt\RG, K_{\wt \RG})$, and let $f'\in \sC_c^\infty(\RG'(\BA))$ be a transfer of $f$.  Then
$$
\Int(f)=-\partial\BJ(f')-\BJ(f'_{\corr}),
$$
where $f'_{\corr}\in \sC_c^\infty(\RG'(\BA))$ is a correction function. Furthermore, we may choose $f'$ such that $f'_{\corr}=0.$
 \end{conj}
\begin{remark}A deeper understanding of the local spherical character $\BJ_{\Pi_v}(f'_v,s)$ (or rather, its derivative) in \eqref{char glo2loc}, together with the spectral decomposition of $\BJ(f',s)$, should allow us to deduce Conjecture \ref{conj:sph AGGP} from Conjecture \ref{conj I=dJ}.
We hope to return to this point in the future.
\end{remark}

The comparison can be localized for a (large) class of test functions $f$ and $f'$.
Let $f=\otimes_v f_v$ be a pure tensor such that there is a place $u_0$ of $F$ where $f_{u_0}$ has support in the regular semisimple locus $\wt\RG(F_{u_0})_\rs$. Then the cycles $\wh R(f)\wh z_{K} $ and $\wh z_{K}$ do not meet in the generic fiber  $\Sh_{K_{\wt \RG}}(\wt \RG)$. The arithmetic intersection pairing then localizes as a sum over all places $w$ of $F'$ (note that $F=\BQ$)
$$
\Int(f)=\sum_{w}\Int_w(f).
$$ 
Here for a non-archimedean place $w$, the local intersection pairing $\Int_w(f)$  is defined through the Euler--Poincar\'e characteristic of a derived tensor product on $\CM_{K_{\wt\RG}}(\wt \RG)\otimes_{O_F} O_{F, w}$, cf.\ \cite[4.3.8(iv)]{GS}. 

Similarly, let $f'=\otimes_v f'_v$ be a pure tensor such that there is a place $u_0$ of $F$ where $f'_{u_0}$ has support in the regular semi-simple locus $\RG'(F_{u_0})_\rs$. Then we have a decomposition 
\begin{align*}
\BJ(f',s)=\sum_{\gamma\in [\RG'(F)_\rs]} \Orb(\gamma,f',s),
\end{align*}
where each term is a product of local orbital integrals \eqref{orb U},
\begin{equation*}
\Orb(\gamma,f',s)=\prod_{v}\Orb(\gamma,f'_v,s).
\end{equation*}
The first derivative $\partial \BJ(f')$ then localizes as  a sum over  places $v$ of $F$,
$$
\partial\BJ(f')=\sum_v\partial\BJ_v(f'),
$$
where   the summand $\partial\BJ_v(f')$ takes the derivative of the local orbital integral (cf.\ \eqref{orb del}) at the place $v$,
$$\partial\BJ_v(f')= \sum_{\gamma\in [\RG'(F)_\rs]} \del(\gamma,f'_v)\cdot \prod_{u\neq v} \Orb(\gamma,f'_u).
$$

It is then natural to expect a place-by-place comparison between $\Int_v(f)=\sum_{w|v}\Int_w(f)$ and $\partial\BJ_v(f')$. If a place $v_0$ of $F$ splits into two places $w_0,\ov w_0$ of $F'$ (and under the above regularity condition on the support of $f$ and of $f'$), we have \cite[Thm. 1.3]{RSZ3}
$$
\Int_{w_0}(f)=\Int_{\ov w_0}(f)=\partial \BJ_{v_0}(f')=0.
$$
For a place $w_0$ of $F'$ above a non-split place $v_0$ of $F=\BQ$, we have a smooth integral model $\CM_{K_{\wt\RG}}(\wt \RG)\otimes_{O_F} O_{F, w_0}$ when $K_{\wt\RG,v_0}$ is a hyperspecial compact open subgroup $\wt \RG(O_{F,v_0})$ (resp.\ a special parahoric subgroup $K^\circ_{\wt\RG,v_0}$) for inert $v_0$ (resp.\ ramified $v_0$).

For an inert $v_0$, the comparison between $\Int_{v_0}(f)$ and $\partial\BJ_{v_0}(f')$ is then reduced to a local conjecture that we will consider in the next section, the {\em arithmetic fundamental lemma} conjecture. 
Let $W^{\flat}[v_0]$ be the pair of nearby Hermitian spaces, i.e., the Hermitian space (with respect to $F'/F$) that is totally negative at archimedean places, and is not equivalent to $W$ at $v_0$. Let $\wt\RG^\flat[v_0]$ be the corresponding group, an inner form of $\wt \RG$.
\begin{thm}\label{thm:A RTF} Let 
$f=\otimes_v f_v$ be a pure tensor such that 
\begin{altenumerate}
  \item $f_{v_0}={\bf 1}_{\wt \RG(O_{F,v_0})}$, and
\item  there is a place $u_0$ of $F$ where $f_{u_0}$ has support in the regular semisimple locus $\wt\RG(F_{v_0})_\rs$.
\end{altenumerate}
Then 
\begin{align}\label{eq A RTF}
\Int_{v_0}(f)=\sum_{g\in \wt\RG^\flat[v_0](F)_\rs}\Int_{v_0}(g)\cdot\prod_{u\neq v_0}\Orb(g,f_u),
\end{align}
where the local intersection number $\Int_{v_0}(g)$ is defined by \eqref{defintprod} in the next section. 
\end{thm}
For $F=\BQ$, this is \cite[Thm. 3.9]{Z12}. The general case is established in (the proof of) \cite[Thm. 8.15]{RSZ3}.

The expansion \eqref{eq A RTF} resembles the geometric side of a usual RTF, and hence we call the expansion (the geometric side of) an {\em arithmetic RTF}.

Finally, for a ramified place $v_0$, the analogous question is reduced to the local {\em arithmetic transfer} conjecture formulated by Rapoport, Smithling and the author in \cite{RSZ1,RSZ2}. We have a result similar to Theorem \ref{thm:A RTF}, \cite[Thm. 8.15]{RSZ3}.

 \subsection{Geometric RTF (over function fields)}
 \label{ss:G RTF}

In the last part of this section, let us briefly recall the strategy to prove the higher Gross--Zagier formula in \S\ref{s:HGZ} over function fields. 

To continue from \S\ref{s:HGZ}, let $f$ be an element in the spherical Hecke algebra $\sH$ (with $\BQ$-coefficient). Let
\begin{align*}
\Int_r(f):=\Big(R(f)*\theta^{\mu}_{*}[\Sht^{\mu}_{T}],\quad \theta^{\mu}_{*}[\Sht^{\mu}_{T}]\Big)_{\Sht'^{\mu}_{G}}
\end{align*}
be the intersection number of the Heegner--Drinfeld cycle with its translation by a Hecke correspondence $R(f)$. Here the right hand side does not depend on $\mu$ but only on the number $r$ of modifications of the Shtukas.

Next, consider  the triple $(\RG',\RH'_1,\RH'_2)$ where 
$\RG'=\RG=\PGL_2$ and $\RH'_1=\RH'_2$ are the diagonal torus $A$ of $\PGL_2$.  Similarly to \ref{J f'}, we define a distribution by a (regularized) integral
\begin{align*}
\BJ(f,s)=\int_{[\RH'_1]}\int_{[\RH_2']}K_{f}(h_1,h_2)\,\lvert h_1h_2\rvert^s \,\eta(h_2)\,dh_1\,dh_2, 
\end{align*}
where, for $h=\diag(a,d)\in A(\BA)$, we write 
$
\big|h\big|=\big|a/d\big|$ and $\eta(h)=\eta_{F'/F}(a/d).$ Let
 $$
\BJ_r(f)=\frac{d^r}{ds^r} \Big|_{s=0} \BJ(f,s).
$$

Then Yun and the author proved in \cite{YZ1} the following {\em key identity}, which we may call a {\em geometric RTF}, in contrast to the arithmetic intersection numbers in the number field case. 
\begin{thm}\label{key id}
Let $f\in\sH$. Then
\begin{equation}\label{key id}
\BI_{r}(f)=(\log q)^{-r}\BJ_{r}(f).
\end{equation}
\end{thm}
In this situation of geometric intersection,  our proof of the key identity \eqref{key id} is entirely global, in the sense that we do {\em not} reduce the identity to the comparison of local orbital integrals. In fact, our proof of \eqref{key id} gives a term-by-term identity of the orbital integrals. This strategy is explained in the forthcoming work of Yun on  the function field analog of the arithmetic fundamental lemma \cite{Yun3}. 

For a more detailed exposition on the geometric construction related to the proof of Theorem \ref{key id}, see Yun's survey \cite{Yun2}.

\section{The arithmetic fundamental lemma conjecture}

We now consider the local version of special cycles on Shimura varieties, i.e., (formal) cycles on Rapoport--Zink formal moduli spaces of $p$-divisible groups. The theorem of Rapoport--Zink on the uniformization of Shimura varieties \cite{RZ} relates the local cycles to the global ones, and this allows us to express the local height of the global cycles (the {\em semi-global} situation in \cite{RSZ3}) to intersection numbers of local cycles, cf.\ Theorem \ref{thm:A RTF}. 

\subsection{The fundamental lemma of Jacquet and Rallis, and a theorem of Yun}
Now let $F$ be a finite extension of $\BQ_p$ for an {\em odd} prime $p$. Let $O_F$ be the ring of integers in $F$, and denote by $q$ the number of elements in the residue field of $O_F$. Let $\breve F$ be the completion of a maximal unramified extension of $F$. Let $F'/F$ be an {\em unramified} quadratic extension.

Recall from \S\ref{sss:JR orb} that there are two equivalence classes of pairs of Hermitian spaces, denoted by $W,W^\flat$ respectively, such that, for $W=(W_{n-1},W_n)$, both $W_{n-1}$ and $W_n$ contain self-dual lattices. We rename the respective groups as $\RG$ and $\RG^\flat$ respectively, and we rewrite the bijection of orbits \eqref{orb bij}:
\begin{align}\label{orb bij1}
\bigl[\RG(F)_{\rm rs}\bigr] \sqcup \bigl[\RG^\flat(F)_{\rm rs}\bigr] \isoarrow  \bigl[\RG'(F)_{\rm rs}\bigr].
\end{align}
Let $\RG(O_F)$ be the hyperspecial compact open subgroup of $\RG(F)$ defined by a self-dual lattice. 
\begin{thm}[Fundamental Lemma (FL)]\label{thm FL}
For a prime $p$ sufficiently large, the characteristic function $\mathbf{1}_{\RG'(O_{F_0})} \in \sC_c^\infty(\RG'(F))$ transfers to the pair of functions $(\mathbf{1}_{\RG(O_F)}, 0) \in \sC_c^\infty(\RG(F)) \times \sC_c^\infty(\RG^\flat(F))$.
\end{thm} 
Jacquet and Rallis conjecture that the same is always true for any {\em odd} $p$ \cite{JR}.  Yun proved the equal characteristic analog of their conjecture for $p> n$; Gordon deduced the $p$-adic case for $p$ large (but unspecified), cf.\ \cite{Yun}.

\subsection{The arithmetic fundamental lemma conjecture}

Next, we let $\CN_n$ be the unitary Rapoport--Zink formal moduli space over $\Spf O_{\breve{F}}$ parameterizing {\em Hermitian supersingular formal $O_{F'}$-modules of signature $(1,n-1)$}, cf.\ \cite{KR-U1,RSZ2}. Let $\CN_{n-1,n}=\CN_{n-1}\times_{\Spf O_{\breve F}}\CN_n$. Then $\CN_{n-1,n}$ admits an action by $\RG^\flat(F)$.

There is a natural closed embedding $\delta \colon \CN_{n-1}\to \CN_n$  (a local analog of the closed embedding \eqref{inc M}). Let 
\begin{equation*}
   \xymatrix{\Delta \colon \CN_{n-1} \xra{\,}   \CN_{n-1,n}}
\end{equation*}
be the graph morphism of $\delta$. We denote by $\Delta_{\CN_{n-1}}$ the image of $\Delta$. For $g\in \RG^\flat(F)_\rs$, we consider the intersection product on $\CN_{n-1, n}$ of $\Delta_{\CN_{n-1}}$ with its translate $g\Delta_{\CN_{n-1}}$, defined through the derived tensor product of the structure sheaves,
\begin{equation}\label{defintprod}
   \Int(g) := \left( \Delta_{\CN_{n-1}}, g\cdot\Delta_{\CN_{n-1}}\right)_{\CN_{n-1, n}} := \chi\left({\CN_{n-1, n}},  \CO_{\Delta_{\CN_{n-1}}}\otimes^\BL\CO_{g\cdot\Delta_{\CN_{n-1}}}\right) . 
\end{equation}

We have defined the derivative of the orbital integral \eqref{orb del}. 
\begin{conj}[Arithmetic Fundamental Lemma (AFL), \cite{Z12}]
\label{U AFL}
Let $\gamma\in\RG'(F)_\rs$ match an element $g\in \RG^\flat(F)_\rs$. Then  
\begin{equation*}\label{introAT}
\omega(\gamma)   \del\left(\gamma,{\bf 1}_{\RG'(O_{F})}\right)=- 2\,\Int(g)\cdot\log q.
\end{equation*}
\end{conj}

 \begin{remark}
\begin{altenumerate}
\item  We may interpret the orbital integrals in terms of ``counting lattices", cf.\ \cite[\S7]{RTZ}. 
\item See \cite[\S4]{RSZ2} for some other equivalent formulations of the AFL conjecture. 
\end{altenumerate}
 \end{remark}

\begin{thm} 
\begin{altenumerate}
  \item 
The AFL conjecture \ref{U AFL} holds when $n=2, 3$. 
  \item When $p\geq\frac{n}{2}+1$, the AFL conjecture \ref{U AFL} holds  for {\em minuscule} elements $g\in \RG^\flat(F)$ in the sense of \cite{RTZ}.
\end{altenumerate}
\end{thm}
Part {\rm (i)} was proved in \cite{Z12}; a simplified proof when $p\geq 5$ is given by Mihatsch in \cite{M-AFL}. Part  {\rm (ii)} was proved by Rapoport, Terstiege, and the author in \cite{RTZ}; a simplified proof is given by Li and Zhu in \cite{LZ}.
 Mihatsch in \cite{M-Th} proved more cases of the AFL for arbitrary $n$ but under restrictive conditions on $g$.
 \begin{remark}
\begin{altenumerate}
  \item Yun has announced a proof of the function field analog of the AFL conjecture \cite{Yun2,Yun3}.
\item  Let $F'/F$ be a {\em ramified} quadratic extension of $p$-adic fields. In \cite{RSZ1,RSZ2}, Rapoport, Smithling, and the author propose an {\em arithmetic transfer} (AT)  conjecture. This conjecture can be viewed as the counterpart of the existence of transfer (cf.\ Theorem \ref{thm ST}) in the arithmetic context over a $p$-adic field.  We proved the conjecture for $n=2,3$.
 \item The analogous question on archimedean local fields remains a challenge, involving Green currents in the complex geometric setting and relative orbital integrals on real Lie groups.  
\end{altenumerate}
 \end{remark}


\begin{thebibliography}{99}


\bibitem{AGRS} {A. Aizenbud, D. Gourevitch, S. Rallis, G. Schiffmann, \textit{Multiplicity one theorems}, Ann. of Math. (2) 172 (2010), no. 2, 1407--1434. }




\bibitem{BP15a}{R. 
Beuzart-Plessis, \textit{
Endoscopie et conjecture raffin\'ee de Gan-Gross-Prasad pour les groupes unitaires}, Compositio Mathematica 151 (2015) 1309--1371.}

\bibitem{BP15b}{\bysame, \textit{A local trace formula for the Gan-Gross-Prasad conjecture for unitary groups: the archimedean case}, arXiv:1506.01452}

\bibitem{BP16}{\bysame, \textit{Comparison of local spherical characters and the Ichino-Ikeda conjecture for unitary groups}, arXiv:1602.06538}




\bibitem{CZ}{P.-H. Chaudouard, M. Zydor, {\it Le transfert singulier pour la formule des traces de Jacquet-Rallis}, arXiv:1611.09656.}



 
\bibitem{De} {P. Deligne, \textit{Travaux de Shimura}. In \textit{S\'eminaire Bourbaki, 23\`eme ann\'ee (1970/71)}, Exp. No. 389, Lecture Notes in Math. \textbf{244}, Springer-Verlag, Berlin, 1971, pp.\ 123--165.}


\bibitem{FJ}{S. Friedberg, H. Jacquet,
\textit{Linear periods.} J. Reine Angew. Math. 443 (1993), 91--139.}



\bibitem{Gan13}{W. T. Gan, {\it Recent progress on the Gross-Prasad conjecture (a survey talk given at the annual meeting of VIASM, July 2013)}, Acta Math Vietnam 39 (2014), no. 1, 11--33.}


\bibitem{GGP}{W. T. Gan, B. Gross, D. Prasad, \textit{Symplectic local root numbers, central critical
$L$-values, and restriction problems in the representation theory of
classical groups},  Ast\'erisque \textbf{346} (2012), 1--109.}



\bibitem{GJR1}{D. Ginzburg, D. Jiang, S. Rallis, \textit{ On the nonvanishing
of the central value of the Rankin-Selberg L-functions.} J. Amer.
Math. Soc. 17 (2004), no. 3, 679--722}


\bibitem{GJR2}{\bysame, \textit{Models for certain residual representations of unitary groups.} Automorphic forms and L-functions I. Global aspects, 125--146, 
Contemp. Math., 488, Amer. Math. Soc., Providence, RI, 2009. }

\bibitem{GS}{H. Gillet and C. Soul\'e, \textit{Arithmetic intersection theory},  Inst. Hautes Etudes Sci. Publ. Math. \textbf{72} (1990), 93--174. } 


\bibitem{G97}{
B. H. Gross, \textit{ On the motive of a reductive group.}
Invent. Math.  130  (1997),  no. 2, 287--313.}

%

\bibitem{GP1}{B. H. Gross, D. Prasad, \textit{On the decomposition of a
representation of ${\rm SO}\sb n$ when restricted to ${\rm SO}\sb
{n-1}$}.  Canad. J. Math.  44  (1992),  no. 5, 974--1002.}

\bibitem{GP2}{\bysame, \textit{On irreducible representations of ${\rm SO}\sb {2n+1}\times{\rm
SO}\sb {2m}$}.  Canad. J. Math.  46  (1994),  no. 5, 930--950.}

\bibitem{GZ}{B. H. Gross, D. Zagier: \textit{Heegner points
and derivatives of $L$-series.} Invent. Math. 84 (1986), no. 2,
225--320. }


\bibitem{HLR}{G. Harder, R. P. Langlands, M. Rapoport, \textit{Algebraische Zyklen auf Hilbert-Blumenthal-Fl\"achen}. J. Reine Angew. Math. 366 (1986), 53--120.}


\bibitem{HN}{
R. N. Harris,  \textit{The refined Gross-Prasad conjecture for unitary groups.} Int. Math. Res. Not. IMRN 2014, no. 2, 303--389.}


\bibitem{Ich}{
A. Ichino, \textit{ Trilinear forms and the central values of triple product
L-functions.}  Duke Math. J. Volume 145, Number 2 (2008), 281--307.}

\bibitem{II}{
A. Ichino, T. Ikeda, \textit{ On the periods of automorphic
forms on special orthogonal groups and the Gross-Prasad
conjecture,} Geom. Funct. Anal. 19 (2010), no. 5, 1378--1425.}



\bibitem{J05}{H. Jacquet, \textit{A guide to the relative trace formula.} In Automorphic representations, L-functions and applications: progress and prospects, 257--272,  Ohio State Univ. Math. Res. Inst. Publ., 11, de Gruyter, Berlin, 2005 }

\bibitem{JR}{H. Jacquet, S. Rallis,\textit{ On the Gross-Prasad conjecture for unitary
groups}. in {\em On certain L-functions}, 205-264, Clay Math. Proc., 13, Amer. Math. Soc., Providence, RI, 2011.}
 
 
\bibitem{Ku02}{S. Kudla, \textit{Derivatives of Eisenstein series and arithmetic geometry.} Proceedings of the International Congress of Mathematicians, Vol. II (Beijing, 2002), 173--183, Higher Ed. Press, Beijing, 2002.}

\bibitem{KR-U1}{S. Kudla, M. Rapoport, \textit{Special cycles on unitary Shimura varieties I. Unramified local theory}, Invent. Math. \textbf{184} (2011), no. 3, 629--682.}

\bibitem{KR-U2}{\bysame, \textit{Special cycles on unitary Shimura varieties II: Global theory},  J. Reine Angew. Math. \textbf{697} (2014), 91--157.}



\bibitem{KRY}{S. Kudla, M. Rapoport, T. Yang, \textit{Modular forms and special cycles on Shimura curves.} Annals of Mathematics Studies, 161. Princeton University Press, Princeton, NJ, 2006. x+373 pp.}

\bibitem{Ku-ens} {K. K\"unnemann, \textit{Height pairings for algebraic cycles on abelian varieties}, Ann. Sci. ENS \textbf{34} (2001), 503--523.}



 \bibitem{La06}{E. Lapid, {\it The relative trace formula and its applications}, Automorphic Forms and Automorphic L-Functions (Kyoto, 2005), Surikaisekikenkyusho Kokyuroku No. 1468 (2006), 76--87.}

 \bibitem{La10}{\bysame, {\it Some applications of the trace formula and the relative trace formula.} Proceedings of the International Congress of Mathematicians. Volume III, 1262--1280, Hindustan Book Agency, New Delhi, 2010}

\bibitem{LZ}{C. Li, Y. Zhu, \textit{Remarks on the arithmetic fundamental lemma}, to appear in Algebra $\&$ Number Theory.}


\bibitem{Liu1}{Y. Liu, \textit{Relative trace formulae toward Bessel and Fourier-Jacobi periods of unitary groups.} Manuscripta Mathematica, 145 (2014) 1--69.}

\bibitem{Liu2}{Y. Liu, \textit{Refined Gan-Gross-Prasad conjecture for Bessel periods.} Journal fur die reine und angewandte Mathematik, 717 (2016) 133--194.}



\bibitem{M-AFL}{A. Mihatsch, \textit{On the arithmetic fundamental lemma conjecture through Lie algebras}, Math. Z. \textbf{287} (2017), no. 1--2, 181--197.}

\bibitem{M-Th}{A. Mihatsch, \textit{Relative unitary RZ-spaces and the arithmetic fundamental lemma}, Ph.D.\ thesis, Bonn, 2016 \href{https://arxiv.org/abs/1611.06520}{\texttt{arXiv:1611.06520 [math.AG]}}.}


\bibitem{MW}{C. M\oe glin, J.L. Waldspurger, {\it La conjecture locale de Gross-Prasad pour les groupes sp\'eciaux orthogonaux: le cas g\'en\'eral}. Sur les conjectures de Gross et Prasad. II. Ast\'erisque No. 347 (2012), 167--216.}


\bibitem{Of09}{O. Offen, {\it Unitary Periods and Jacquet's relative trace formula.} Automorphic forms and L-functions, Proceedings of a workshop in honor of Steve Gelbart on the occasion of his 60th birthday, Contemporary Mathematics 488, AMS and BIU (2009) 148--183.}

\bibitem{RSZ1}{M. Rapoport, B. Smithling, W. Zhang, \textit{On the arithmetic transfer conjecture for exotic smooth formal moduli spaces}, Duke Math. J. \textbf{166} (2017), no. 12, 2183--2336.} 

\bibitem{RSZ2}{\bysame, \textit{Regular formal moduli spaces and arithmetic transfer conjectures}, to appear in Math. Ann.}


\bibitem{RSZ3}{\bysame, \textit{Arithmetic diagonal cycles on unitary Shimura varieties}, preprint, arXiv:1710.06962.}


\bibitem{RTZ}{M. Rapoport, U. Terstiege, W. Zhang, \textit{On the arithmetic fundamental lemma in the minuscule case}, Compos. Math. \textbf{149} (2013), no. 10, 1631--1666.}


\bibitem{RV}{M. Rapoport, E. Viehmann,  \textit{Towards a theory of local Shimura varieties,}
Munster J. of Math. 7 (2014), 273--326.}

\bibitem{RZ}{M. Rapoport, Th. Zink, \textit{Period spaces for $p$-divisible groups}, Annals of Mathematics Studies, vol. \textbf{141}, Princeton University Press, Princeton, NJ, 1996.}


\bibitem{Sa08}{Y. Sakellaridis, \textit{On the unramified spectrum of spherical varieties over $p$-adic fields.} Compositio Mathematica 144 (2008), no. 4, 978--1016}

\bibitem{SV}{Y. Sakellaridis, A. Venkatesh, \textit{Periods and harmonic analysis on spherical varieties.} To appear in Ast\'erisque, 296pp, arXiv:1203.0039.}



\bibitem{SZ}{B. Sun, C. Zhu, \textit{Multiplicity one theorems: the Archimedean case}, Ann. of Math. (2) 175 (2012), no. 1, 23--44.} 



\bibitem{W85}{J. Waldspurger, \textit{Sur les valeurs de certaines fonctions
L automorphes en leur centre de sym\'etrie.}  Compositio Math. 54
(1985), no. 2, 173--242.}


\bibitem{W12}{\bysame, \textit{La conjecture locale de Gross-Prasad pour les repr\'esentations temp\'er\'ees des groupes sp\'eciaux orthogonaux.} Sur les conjectures de Gross et Prasad. II. Ast\'erisque No. 347 (2012), 103--165. }


\bibitem{Xue}{H. Xue, \textit{On the global Gan--Gross--Prasad conjecture for unitary groups: approximating smooth transfer of Jacquet--Rallis}, To appear in J. Reine Angew. Math.}


\bibitem{YZZ1}{X. Yuan, S.-W. Zhang, W. Zhang, \textit{The Gross--Zagier formula on Shimura curves}, Annals of Mathematics Studies, vol. \textbf{184}, Princeton University Press, Princeton, NJ, 2013.}

\bibitem{YZZ2}{\bysame, \textit{Triple product $L$-series and Gross--Kudla--Schoen cycles}, in preparation.}

\bibitem{Yun}{Z. Yun, \textit{The fundamental lemma of Jacquet--Rallis in positive
characteristics}, with an appendix by Julia Gordon,  Duke Math. J. \textbf{156} (2011), no. 2, 167--228.}
\bibitem{Yun2}{\bysame, \textit{\it Hitchin type moduli stacks in automorphic representation theory}, ICM report.}



\bibitem{Yun3}{\bysame, \textit{\it The arithmetic fundamental lemma over function fields}, in preparation.}


\bibitem{YZ1}{Z. Yun, W. Zhang, \textit{Shtukas and the Taylor expansion of L-functions}, Ann. of Math., Volume 186 (2017), no. 3, pages 767--911. }


\bibitem{YZ2}{\bysame, \textit{Shtukas and the Taylor expansion of L-functions (II)}, preprint. }



\bibitem{ZSW10}{S. Zhang. \textit{Linear forms,  algebraic cycles, and derivatives of L-series}, preprint, 2010.}




\bibitem{Z09}{W. Zhang, \textit{Relative trace formula and arithmetic Gross--Prasad conjecture}, unpublished manuscript, 2009.}

\bibitem{Z12}{\bysame, \textit{On arithmetic fundamental lemmas}, Invent. Math. \textbf{188} (2012), no. 1, 197--252.}

\bibitem{Z12b}{\bysame, \textit{Gross-Zagier formula and arithmetic fundamental lemma.} In \textit{Fifth International Congress of Chinese Mathematicians. Part 1, 2}, AMS/IP Stud. Adv. Math. \textbf{51}, pt. 1, 2, Amer. Math. Soc., Providence, RI, 2012, pp. 447--459.} 


\bibitem{Z14a}{\bysame, \textit{Fourier transform and the global Gan--Gross--Prasad conjecture for unitary groups},  Ann. of Math. (2) \textbf{180} (2014), no. 3, 971--1049.}

\bibitem{Z14b}{\bysame, \textit{Automorphic period and the central value of Rankin-Selberg $L$-function}, J. Amer. Math. Soc. \textbf{27} (2014), 541--612.}

\bibitem{Z17}{\bysame, \textit{A conjectural linear arithmetic fundamental lemma for Lubin--Tate spaces}, unpublished manuscript, 2017}


\bibitem{Zy}{M. Zydor, \textit{Les formules des traces relatives de Jacquet-Rallis grossi\`eres}, 
preprint, 2015, \href{http://arxiv.org/abs/1510.04301}{\texttt{arXiv:1510.04301 [math.NT]}}.}

\end{thebibliography}
\end{document}